\input amstex
\input epsf
\define\labelnumber{0}
\ifx\labelsloaded\relax \else\let\labelsloaded\relax\fi
\catcode`\@=11 
\newwrite\@auxout
\newif\if@fwdref \@fwdreffalse
 
%
\newcount\Ch@nro \Ch@nro=0
\newcount\Sec@nro \Sec@nro=0
\newcount\Tag@nro \Tag@nro=0
\newcount\Thm@nro \Thm@nro=0
\newcount\Fig@nro \Fig@nro=0
\newcount\Table@nro \Table@nro=0
\newif\ift@gsandthms \t@gsandthmsfalse
%
%
\def\theCh{\the\Ch@nro}
\def\theSec{\the\Sec@nro}
\def\theTag{\ch.\the\Tag@nro}
\def\theThm{\ch.\the\Thm@nro}
\def\theFig{\ch.\the\Fig@nro}
\def\theTable{\the\Table@nro}
\def\TagsAndTheorems{\t@gsandthmstrue}
\def\stepbyone#1{\global\expandafter
  \advance\csname#1@nro\endcsname by1\relax}
\def\stepcounter#1#2{\global\expandafter
  \advance\csname#1@nro\endcsname by #2\relax}
\def\setcounter#1#2{\global\expandafter\csname#1@nro\endcsname =#2\relax}
 
\def\Ch@{\setcounter{Sec}{0}\setcounter{Fig}{0}%
  \setcounter{Tag}{0}\setcounter{Thm}{0}%
   \xdef\@currentlabel{\ifnum\Ch@nro=0 \else\theCh\fi}\@currentlabel}
\def\Ch{\stepbyone{Ch}\Ch@}

\def\Sec@{%
  \xdef\@currentlabel{\ifnum\Ch@nro=0 \else\theCh.\fi
     \ifnum\Sec@nro=0 \else\theSec\fi}%
  \@currentlabel}
\def\Sec{\stepbyone{Sec}\Sec@}
 
\def\LeftTagForm{(}
\def\RightTagForm{)}
\def\Tag@{\xdef\@currentlabel{\LeftTagForm
    \ifnum\Ch@nro=0 \else\theCh.\fi
    \theTag\RightTagForm}%
  \@currentlabel}
\def\Tag{\stepbyone{Tag}%
    \ift@gsandthms\else\stepbyone{Thm}\fi\Tag@}
 
\def\Thm@{\xdef\@currentlabel{\ifnum\Ch@nro=0 \else\theCh.\fi
    \ifnum\Sec@nro=0 \else\theSec.\fi\theThm}%
  \@currentlabel}
\def\Thm{\stepbyone{Thm}%
    \ift@gsandthms\else\stepbyone{Tag}\fi\Thm@}
 
\def\Fig@{\xdef\@currentlabel{\theFig}\@currentlabel}
\def\Fig{\gdef\CorrectCounter{\stepcounter{Fig}{-1}}%
  \stepbyone{Fig}%
  \Fig@}
 
\def\Table@{\xdef\@currentlabel{\theTable}\@currentlabel}
\def\Table{\gdef\CorrectCounter{\stepcounter{Table}{-1}}%
  \stepbyone{Table}\Table@}
%
 
%
%
 
\def\thepage{\number\pageno}
 
\def\refer#1{\@ifundefined{r@#1}{{ ??}\global\@fwdreftrue\@warning
   {Reference `#1' on page \thepage \space
    undefined}}{\edef\@tempa{\@nameuse{r@#1}}\expandafter
    \@car\@tempa \@nil\null}}
 
\def\referpage#1{\@ifundefined{r@#1}{{ ??}\global\@fwdreftrue\@warning
   {Reference `#1' on page \thepage \space
    undefined}}{\edef\@tempa{\@nameuse{r@#1}}\expandafter
    \@cdr\@tempa\@nil\null}}
 
\newif\if@filesw \@fileswtrue
 
\def\marginal#1{\strut\setbox0=%
      \vtop{\hsize=15mm
            \sevenrm
            \textfont0=\scriptfont0 \textfont1=\scriptfont1
            \textfont2=\scriptfont2 \textfont3=\scriptfont3
            \parindent=0pt\baselineskip=8pt
             \raggedright\rightskip=1em plus 2em\hfuzz=.5em\tolerance=9000
             \overfullrule=0pt
             #1}%
      \dimen0=\ht0
      \advance\dimen0 by \dp\strutbox
      \ht0=0pt\dp0=\dimen0
      \vadjust{\kern-\dimen0\moveleft\wd0\box0}%
      \ignorespaces}       

\def\maybemarginal#1{\ifnum\labelnumber=1{\marginal{#1}\igs}\else{}\fi}

\def\label#1{\maybemarginal{#1}\igs\@bsphack\if@filesw {\let\thepage\relax
   \xdef\@gtempa{\write\@auxout{\string
      \newlabel{#1}{{\@currentlabel}{\thepage}}}}}\@gtempa
   \fi\@esphack\igs}

\def\rosterlabel#1{\@bsphack\xdef\@currentlabel{%
   \therosteritem{\number\rostercount@}}%
   \if@filesw {\let\thepage\relax
   \xdef\@gtempa{\write\@auxout{\string
      \newlabel{#1}{{\@currentlabel}{\thepage}}}}}\@gtempa
   \fi\@esphack}
 
\def\newlabel#1#2{\@ifundefined{r@#1}{}{\@warning{Label `#1' multiply
   defined}}\global\@namedef{r@#1}{#2}}
 
\def\@currentlabel{} 
 
%
\newif\if@indexsw
\newwrite\@idxout
 
\def\makeindex{\@indexswtrue\immediate\openout\@idxout=\jobname.idx
   \immediate\write\@idxout{\string\relax}}
 
\def\index#1{\@bsphack\if@indexsw {\let\thepage\relax
   \xdef\@gtempa{\write\@idxout{\string
      \indexentry{#1}{\thepage}}}}\@gtempa
   \fi\@esphack}
 
\def\indexentry#1#2{\noindent {#1}{\unskip\nobreak\hfil\penalty50
   \hskip2em\hbox{}\nobreak\hfil#2%
   \parfillskip=0pt \finalhyphendemerits=0 \par}}
 
%
%
\def\make@ref#1{\@ifundefined{b@#1}{{??}\global\@fwdreftrue\@warning
   {Citation `#1' on page \thepage \space
    undefined}}{\@nameuse{b@#1}}}
%
%
\def\citeto#1{\relaxnext@
 \def\nextiii@##1,##2\end@{\immediate\write\@auxout{\string\citation{##1}}%
   [{\bf\make@ref{##1}},##2]}%
 \in@,{#1}\ifin@\def\next{\nextiii@#1\end@}\else
   \def\next{\immediate\write\@auxout{\string\citation{#1}}%
     {\bf\make@ref{#1}}}\fi
     \next}
\def\noOf#1{\def\name{#1}
   \global\advance\c@num by1
   \if@filesw \immediate\write\@auxout
   {\string\bibcite{#1}{\the\c@num}}\fi\no\make@ref{#1}}
\def\keyOf#1#2{\if@filesw
   {\def\protect##1{\string ##1\space}\immediate
   \write\@auxout{\string\bibcite{#1}{#2}}}\fi
   \key{$\lbrack${\bf #2}$\rbrack$}}
 
\def\TAG{\tag"\Tag"}
\let\@@end=\end
\let\@@supereject=\supereject
\def\supereject{\@@supereject}
\def\end{\if@fwdref
  \typeout{Undefined forward references/citations. Rerun.}\fi
  \@@end}
%
%
 
%
%
\newwrite\@unused
\def\typeout#1{{\let\protect\string\immediate\write\@unused{#1}}}
\def\@warning#1{\typeout{Warning: #1.}}
 
\def\@namedef#1{\expandafter\def\csname #1\endcsname}
\def\@nameuse#1{\csname #1\endcsname}
 
\def\@car#1#2\@nil{#1}
\def\@cdr#1#2\@nil{#2}
 
\long\def\@ifundefined#1#2#3{\expandafter\ifx\csname
  #1\endcsname\relax#2\else#3\fi}
 
\def\@ifnextchar#1#2#3{\let\@tempe #1\def\@tempa{#2}\def\@tempb{#3}\futurelet
    \@tempc\@ifnch}
\def\@ifnch{\ifx \@tempc \@sptoken \let\@tempd\@xifnch
      \else \ifx \@tempc \@tempe\let\@tempd\@tempa\else\let\@tempd\@tempb\fi
      \fi \@tempd}
 
\def\?{\let\@sptoken= } \?  
 
\def\?{\@xifnch} \expandafter\def\? {\futurelet\@tempc\@ifnch}
 
%
\newdimen\@savsk
\newcount\@savsf
 
\def\@bsphack{\@savsk\lastskip
    \ifhmode\@savsf\spacefactor\fi}
 
\def\@esphack{\relax\ifhmode\spacefactor\@savsf
     {}\ifdim \@savsk >\z@ \ignorespaces\fi\fi}
 
%

%
\newcount\c@num \c@num=0
 
\def\bibcite#1#2{\@ifundefined{r@#1}{}{\@warning{Referencelabel `#1' multiply
   defined}}\global\@namedef{b@#1}{#2}}
 
\def\@gobble#1{}
\let\citation\@gobble
 
%
%
%
\def\ins@#1{\relaxnext@
 \smallcaptionwidth@\captionwidth@\gdef\thespace@{#1}%
 \def\next@{\ifx\next\space@\def\next. {\futurelet\next\nextii@}\else
  \def\next.{\futurelet\next\nextii@}\fi\next.}%
 \def\nextii@{\ifx\next\caption\def\next\caption{\futurelet\next\nextiii@}%
  \else\let\next\nextiv@\fi\next}%
 \def\nextiv@{\vnonvmode@
  {\ifmid@\let\next\midinsert\else\let\next\topinsert\fi
  \next\vbox to\thespace@{}\endinsert}
  {\ifmid@\nonvmodeerr@\midspace\else\nonvmodeerr@\topspace\fi}}%
 \def\nextiii@{\ifx\next\captionwidth\let\next\nextv@
  \else\let\next\nextvi@\fi\next}%
 \def\nextv@\captionwidth##1##2{\smallcaptionwidth@##1\relax\nextvi@{##2}}%
 \def\nextvi@##1{\def\thecaption@{##1}%
  \def\next@{\ifx\next\space@\def\next. {\futurelet\next\nextvii@}\else
   \def\next.{\futurelet\next\nextvii@}\fi\next.}%
  \futurelet\next\next@}%
 \def\nextvii@{\vnonvmode@
  {\ifmid@\let\next\midinsert\else
  \let\next\topinsert\fi\next\vbox to\thespace@{}\nobreak\smallskip
  \setbox\z@\hbox{\noindent\ignorespaces\thecaption@\unskip}%
  \ifdim\wd\z@>\smallcaptionwidth@\centerline{\vbox{\hsize\smallcaptionwidth@
  \noindent\ignorespaces\CorrectCounter\thecaption@
  \gdef\CorrectCounter{}\unskip}}%
  \else\centerline{\box\z@}\fi\endinsert}
  {\ifmid@\nonvmodeerr@\midspace
  \else\nonvmodeerr@\topspace\fi}}%
 \futurelet\next\next@}
%
%
\openin1 \jobname.aux \ifeof1 \typeout{No file \jobname.aux.}\closein1
    \else\closein1 \relax\input \jobname.aux \fi
\immediate\openout\@auxout=\jobname.aux
\immediate\write\@auxout{\string\relax}
\catcode`\@=\active

\documentstyle{amsppt}

\catcode`\@=11
 \def\logo@{}
 \def\captionfont@{\eightpoint\rm}
\catcode`\@=13

 1
 1

\magnification=1200
\def\section#1{\setcounter{Tag}{0}\setcounter{Thm}{0}\setcounter{Fig}{0}\igs
   \goodbreak\flushpar{\bf #1} \nobreak
  \bigskip \nobreak\flushpar\ignorespaces}
\def\endsection{\bigskip}
\def\subsection#1{\medskip\flushpar{\smc #1}\smallskip
                  \flushpar\ignorespaces}

\def\today{\number\day
\space\ifcase\month\or
  January\or February\or March\or April\or May\or June\or
  July\or August\or September\or October\or November\or December\fi
  \space\number\year}

\def\wt#1{\widetilde{#1}}

\def\PP{{\Bbb P}}

\def\qed{{\hfill $\square$\nl}}
\def\es{\varnothing}

\def\ii{\itemitem}
\def\ci{\citeto}

\def\q{\quad}
\def\qq{\qquad}

\def\la{\langle}
\def\ra{\rangle}

\def\s{\sigma}

\def\g{\gamma}
\def\b{\beta}
\def\a{\alpha}
\def\t{\theta}
\def\La{\Lambda}

\def\sT{{\Cal T}}
\define\sF{{\Cal F}}
\def\sR{{\Cal R}}
\def\sB{{\Cal B}}
\def\sE{{\Cal E}}
\def\sEf{\sE_{\text{\rm o}}}
\def\sEo{\sE_{\text{\rm o}}}
\def\sEm{\sE_{\text{\rm m}}}
\def\sEc{\sE_{\text{\rm c}}}
\def\sEcu{\sE_{\text{\rm cut}}}
\def\sSR{{\Cal S}}
\def\sSRm{\sSR_{\text{\rm B}}}
\def\sSRf{\sSR_{\text{\rm o}}}
\def\sSRc{\sSR_{\text{\rm c}}}
\def\pd{\partial}
\def\bop={{\bold p}}

\def\pga{{p_G}}

\def\sC{{\Cal C}}

\def\Om{\Omega}
\def\om{\omega}
\comment
\def\marginal#1{\strut\setbox0=%
      \vtop{\hsize=15mm
            \sevenrm
            \textfont0=\scriptfont0 \textfont1=\scriptfont1
            \textfont2=\scriptfont2 \textfont3=\scriptfont3
            \parindent=0pt\baselineskip=8pt
             \raggedright\rightskip=1em plus 2em\hfuzz=.5em\tolerance=9000
             \overfullrule=0pt
             #1}%
      \dimen0=\ht0
      \advance\dimen0 by \dp\strutbox
      \ht0=0pt\dp0=\dimen0
      \vadjust{\kern-\dimen0\moveleft\wd0\box0}%
      \ignorespaces}       
\endcomment

\def\R{{\Bbb R}}
\def\N{{\Bbb N}}

\def\pep{P_{E,p}}
\def\pgf{probability generating function}
\def\igs{\ignorespaces}
\redefine\qed{\hfill$\square$}
\def\oo{\infty}
\def\capt#1{\baselineskip=10pt\stepcounter{Fig}{-1}\botcaption{\baselineskip=10pt\eightpoint{\bf Figure \Fig.}
  #1}\endcaption\par}
\def\mletter#1#2#3{\hskip#2cm\lower#3cm\rlap{$#1$}\hskip-#2cm}
\def\lastletter#1#2#3{\hskip#2cm\lower#3cm\rlap{$#1$}\hskip-#2cm\vskip-#3cm}

\def\figure#1\par{\parindent=0pt
  \vbox{\baselineskip=0pt \lineskip=0pt
  \line{\hfil}
  #1}}
\def\rc{random-cluster}
\def\Zlpqxs{Z^{\xi,\sim}_{\Lpq}}
\def\Lpq{\Lambda,p,q}
\def\sm{\setminus}
\def\resp{respectively}
\def\lra{\leftrightarrow}

\def\nlra{\nleftrightarrow}
\def\bin{\text{\rm bin}}
\def\pc{p_{\text{\rm c}}}
\def\pg{p_{\text{\rm b}}}
\def\yellow{yellow}
\def\green{black}

\def\pik{p_i(k;i_1,\dots,i_k)}
\def\ik{i_1,\dots,i_k}
\def\qij{q_{ij}}
\def\qijl{q_{i_l j_l}}
\def\Phik{\Phi(k;i_1,j_1,\dots,i_k,j_k)}
\def\ijk{(i_1,j_1),\dots,(i_k,j_k)}
\def\pijk{p_{ij}(k;i_1,j_1,\dots,i_k,j_k)}
\def\sI{{\Cal I}}
\def\sJ{{\Cal J}}
\def\set#1{\{#1\}}
\def\De{\Delta}
\def\eps{\epsilon}
\def\ol#1{\overline{#1}}
\def\gest{\ge_{\text{\rm st}}}
\def\lest{\le_{\text{\rm st}}}
\def\ins{\text{\rm ins}}
\def\bigmid{\,\big|\,}
\def\cmax{C_{\text{\rm max}}}
\def\crw{\sR_{W}^{\text{\rm c}}}
\define\edc{E_{\Delta,C}}
\define\K{{\ol K}}

\define\xrho{\pi}
\define\xsR{\Pi}

\topmatter
\title
Branching Processes, and Random-Cluster Measures on Trees
\endtitle
\rightheadtext{Branching Processes, and Random-Cluster Measures on Trees}

\author
Geoffrey Grimmett, Svante Janson
\endauthor

\address 
Geoffrey Grimmett, Statistical Laboratory, Centre for
Mathematical Sciences, University of Cambridge,
Wilberforce Road, Cambridge
CB3 0WB, UK
\endaddress
\email g.r.grimmett{\@}statslab.cam.ac.uk
\endemail
\address 
Svante Janson, Department of Mathematics, Uppsala University, PO Box
480, S-751 06 Uppsala, Sweden
\endaddress
\email svante.janson{\@}math.uu.se \endemail

\abstract
Random-cluster measures on infinite regular trees
are studied in conjunction with a general
type of `boundary condition', namely an equivalence
relation on the set of infinite paths of the tree.
The uniqueness and non-uniqueness of \rc\ measures
are explored for certain classes of equivalence 
relations. In proving uniqueness, the following problem
concerning branching processes is encountered
and answered.
Consider bond percolation on the family-tree $T$ of
a branching process. What is the probability that every
infinite path of $T$, beginning at its root, 
contains some 
vertex which is itself the root of an infinite open 
sub-tree? 
\endabstract



\endtopmatter\footnote""{This version was prepared on \today.}

\document
\def\ch{1}
\section{\ch. Introduction and summary}
The random-cluster model may be
viewed as a unification of percolation and the 
Ising/Potts models for a ferromagnet. It was described
by Fortuin and Kasteleyn around 1970 in a series of 
influential papers, and has provided one of the principal 
methods for studying the mathematics of ferromagnetism.
 See
[\ci{G02}, \ci{G04}] for detailed accounts
of and
bibliographies associated with the model. When 
the underlying graph $G$ is finite, the 
corresponding \rc\  measure is given in a closed form.
When the graph is infinite, one proceeds either by taking
weak limits of measures on finite subgraphs $H$
as $H\uparrow G$, or by concentrating on a class
of measures whose conditional measures, given
the configuration off a finite subgraph $H$, satisfy
the appropriate  `DLR/Gibbs specification'
(see [\ci{G93}]). Much (but not all) is known
about the relationship between these two approaches when $G$ is a
finite-dimensional lattice. 
The primary purpose of the present
paper is to study the corresponding problem
when $G$ is an infinite regular tree, thus continuing 
a project initiated in [\ci{Hag96}]. 

A random-cluster measure on a finite tree is simply a product measure
--- it is the circuits 
in a graph which cause dependence between the
states of different edges, and, when there are no 
circuits,
there is no dependence. Circuits
may, however, be introduced into trees through a consideration of
boundary conditions, and there
 lies the principal direction
of this paper.  Let $T$ be an infinite labelled tree
with root 0, and let $\sR$ be the set of all infinite
(self-avoiding) paths of $T$ beginning at 0,
termed {\it rays\/}. We may think
of a boundary condition as being
an equivalence relation $\sim$
on $\sR$, the `physical' meaning 
being that two rays $\rho$, $\rho'$ are considered to
be `connected at infinity' whenever $\rho\sim\rho'$.
Such connections affect the counts of 
connected components of random subgraphs, thereby 
contributing
to the random-cluster measures defined on $T$.
The two extremal boundary conditions are usually termed
`free' (meaning that there exist no connections at 
infinity) and
`wired' (meaning that all rays are equivalent),
respectively; these notions agree with those in use
for lattices. The wired boundary condition on $T$
is that studied in [\ci{Hag96}, \ci{Jon99}].

Our study of \rc\ measures will be pursued in
Sections 4--6 in the context of the infinite
$m$-ary tree $T_m'$, where $m\in\{2,3,\dots\}$.
Let $\sim$ be an equivalence relation on the set
$\sR$ of rays of $T_m'$. In Section 4 is presented the DLR/Gibbs specification of a so-called
($\sim$){\it \rc\ measure\/} on $T$. When studying 
\rc\ measures
which arise through limits of finite-volume measures, it
turns out to be natural to restrict oneself to equivalence classes
which are `closed' when viewed as subsets of $\sR^2$.
Thus we are led to consider the topological properties
of equivalence relations, and this we do in Section 5.

A \rc\ 
measure $\phi_{G,p,q}$ on a graph $G$
has two parameters, namely an 
edge-weight $p\in[0,1]$ and a cluster-weight 
$q\in(0,\oo)$. It is an important and useful 
property of \rc\ measures with $q\ge 1$ that they
satisfy the FKG inequality, and for this reason we 
confine ourselves here to this case.
The measure $\phi_{G,p,q}$ increases (in the sense
of stochastic ordering) as $p$ increases. When $G$
is infinite, one knows in the
case of lattices (see [\ci{G93}])
that there exists a unique \rc\ measure 
with parameters $p$
and $q$ ($\ge 1$) whenever $p$
is either sufficiently small or sufficiently
large, and it is an important open problem to
determine the uniqueness region exactly.
The case of small $p$ was answered for  
$T_m'$ in
[\ci{Hag96}], where uniqueness was proved for
all $p< p_{m,q}$ where $p_{m,q}$
is given by an explicit formula.
It was proved moreover that there
exists an interval of values of $p$ 
of the form $[p_{m,q},
p'_{m,q}]$,  non-empty when $q >2$, such
that there is non-uniqueness of wired
\rc\ measures for $p$ lying in
this interval. It was conjectured in [\ci{Hag96}]
that uniqueness of wired measures
is valid when $p > p'_{m,q}$, and such uniqueness was
proved in [\ci{Jon99}] for
sufficiently large $p$.
In Section 6, we extend the work
of [\ci{Jon99}] to prove the existence of
$p''_{m,q} \in (p'_{m,q},1)$ such that
uniqueness is valid for $p \ge p''_{m,q}$
in the more general context of
a certain sub-class of equivalence relations
termed `open' relations.

In proving the above uniqueness, we make use of
a result from branching processes which may have other applications also. Let $T$ be the family-tree of a Galton--Watson branching process with a
single progenitor 
0, and we assume for simplicity that
every family-size is at least 1 and that the mean family
size exceeds 1. On $T$ we construct a bond percolation
process with given edge-density $p$. A vertex $v$ of
$T$ is coloured {\it blue\/} if it is the root
of an infinite open sub-tree of $T$. The progenitor
0 is coloured {\it\green\/} if every infinite path of
$T$ starting at 0 contains some blue vertex.
We shall see in Section 2 how 
to calculate the probability $\g=
\PP(\text{0 is \green})$  
in terms of the family-size
\pgf\ $G$ of the branching process.

Of special relevance to our study of \rc\ measures
is the problem of finding a necessary and sufficient
condition on $p$ such that $\g = 1$. We shall see
in Theorem 2.2 that $\g=1$ if and only if $p\ge \pga$,
where $\pga$ is given uniquely by the equation
$G'(1-\pga \t(\pga))=1$, and 
$\t(p)$ is the survival probability
of the open sub-tree with root 0. [We consider here
the `quenched' probability measure which
`averages' over realizations of $T$ as well as
over the open edge-set of $T$.] Although we
obtain such results in the context of a general 
branching process, in our application to the \rc\ model,
we shall consider (as in [\ci{Jon99}])
only the deterministic case in which every individual 
has exactly $m$ children; this is the
case
with $G(x)=x^m$.

The present work is related to the analysis of 
the \rc\  model on
the complete graph performed in [\ci{BGJ}] and
continued in  [\ci{LL04}], the common concept being 
that of a `mean-field model'. 
A mean-field theory of statistical mechanics arises
either through removing the finite-dimensional
spatial aspect of the system, or by
considering a model which is in some sense 
`infinite-dimensional'. In seeking rigorous
theory, mathematicians often consider the correct
setting for a mean-field model to be either the complete graph or a
tree. In the case of percolation, for example, the corresponding
models are the so-called  
Erd\H os--R\'enyi random graph
(see [\ci{Bol}, \ci{JLR}]) and the binomial branching process
(see [\ci{G99}, Chapter 10]).   
Paper [\ci{BGJ}] contains the
theory of the \rc\ measure on the complete graph with $n$ vertices,
where $q\in(0,\oo)$ and $p=\lambda/n$, in the limit as $n\to\oo$.
The present paper (and the earlier [\ci{Hag96},
\ci{Jon99}]) 
is devoted to the
case of trees.
\endsection

\def\ch{2}
\section{\ch. Branching processes}
We pose and answer a natural question concerning branching
processes. This has an application for the uniqueness of \rc\ measures
on trees, and it may well have further applications
in other areas of probability theory.
It may be viewed as an extension of
a (sub-)result of [\ci{Jon99}].

We consider a (Galton--Watson) branching process with family-size \pgf\ $G$ satisfying
$$
G(0)=0, \q 1<G'(1)<\oo.
\TAG
$$\label{2.ass}
In other words, the number of offspring
of any individual is non-zero, and the mean family-size is strictly greater than 1 (we shall return 
after Corollary
2.3 to the situation in which one dispenses with
the assumption $G(0)=0$). 
The family-tree $T$ of the process 
is an infinite tree with
labelled vertex-set $V$ and a single progenitor called its {\it origin\/}
and labelled 0. We write $\PP$ for the probability
measure which governs the branching process.
For general accounts of the theory of branching processes,
see [\ci{AN72}, \ci{Ha02}, \ci{Jag75}].

We turn $T$ into a directed tree by directing every edge away from 0.
Let $x$ be a vertex. An {\it $x$-ray\/} is defined to be
an infinite directed path of $T$ with (unique) endvertex $x$.
We denote by $\sR_x$ the set of all $x$-rays of $T$, and we abbreviate
$\sR_0$ to $\sR$. We shall use the term {\it ray\/} to mean a member of 
some $\sR_x$. The edge of $T$ joining vertices
$x$ and $y$ is denoted $\langle x,y\rangle$ when
undirected, and $[x,y\rangle$ when directed from $x$ to $y$.

We introduce next a second level of randomness through a consideration
of bond percolation on $T$ (see [\ci{G99}] for a general account of percolation).
Suppose for the moment that $T$ is given, which is to say that $T=(V,E)$
is a given labelled directed rooted infinite tree as above.
Let $0\le p\le 1$. Each edge of $T$ is declared {\it open\/} with
probability $p$, and {\it closed\/} otherwise; the states of
different edges are independent. This amounts to considering
the product measure $\pep$ with density $p$
on the configuration space $\{0,1\}^E$. Let $\wt E$ ($\subseteq E$) be 
the set of open edges, and define the forest $\wt T=(V,\wt E)$.
It is standard that the connected component of $\wt T$ containing 0 is itself
a branching process. A path $\pi$ of $T$ is called {\it open\/} if every edge
in $\pi$ is open.

The vertices of $T$ are 
assigned colours depending
on the sub-trees of which they are roots.
Let $x\in V$. We colour $x$ {\it blue\/} if some $x$-ray of $T$ is
open; we colour $x$ {\it \yellow\/} if $x$ is not blue but every $x$-ray
of $T$ contains some blue vertex; we colour $x$ {\it red\/}
if it is neither blue nor \yellow. Finally, a vertex which is either blue
or \yellow\ is said to be {\it \green\/}. Note that $x$ is \green\
if and only if every $x$-ray contains a blue vertex. We write
$$
\g_{T,x} = \pep(\text{$x$ is \green}),\q
\t_{T,x} = \pep(\text{$x$ is blue}),
\TAG
$$\label{2.defg}
noting that these quantities depend on the tree $T$. We now
average over the measure $\PP$.
Let $k \ge 0$ and let $\sF_k$
be the $\s$-field generated by the first $k$ generations of $T$.
Suppose that $v$ lies in the $k$th generation of $T$. By the
Markov property of branching processes, the quantities
$$
\g_v=\PP(\g_{T,v}\mid \sF_k),\q \t_v=\PP(\t_{T,v}\mid \sF_k)
\TAG
$$\label{2.defg2}
do not depend (almost surely) on the choice of $v$ and $k$, whence
$$
\g_v = \g_0,\q \t_v=\t_0.
\TAG
$$\label{2.defg3}
[Rather than introduce further notation, 
we use $\mu(X)$ to denote the mean of a random
variable $X$ under a
probability measure $\mu$.]
We introduce the abbreviations
$$
\g=\g_0,\q \t=\t_0,
\TAG
$$\label{2.defg4}
and we note the obvious inequality
$$
\g \ge \t.
\TAG
$$\label{2.defg5}
In summary, the root 0 is blue (\resp, red, yellow) with
probability $\t$ (\resp, $1-\g$, $\g-\t$); it is black
with probability $\g$.

The calculation of $\t=\t(p,G)$ is standard, and may be found in any
of many textbooks 
(see, for example, [\ci{GS01}, Thm 5.4.5]).
The extinction probability $\eta=1-\t$ is the smallest positive root
of the equation 
$$
\eta = G(1-p+p\eta),
\TAG
$$\label{2.extcond}
and thus $\t$ is the largest root in $[0,1]$
of the equation
$$
\t = 1 - G(1-p\t).
\TAG
$$\label{2.extcond2}
It follows from \refer{2.extcond}
in the usual manner that
$$
\t>0 \q\text{if and only if}\q 
pG'(1)>1.
\TAG
$$\label{2.extcond3}

Our principal target in this section is to
calculate $\g$. We define $f_p:[p\t,1]
\to [0,\oo)$ by
$$
f_p(\a)=\t+G(\a-p\t),\qq \a\in[p\t,1],
\TAG
$$\label{2.defnf}
and we note some properties of $f_p$ the 
proofs of which will
be given later.

\proclaim{Proposition \ch.1}
\label{2.prop1} 
Let $p\in[0,1)$ and let $G$ satisfy
\refer{2.ass}.
The equation $\a=f_p(\a)$ has a root at
$\a=1$. It has either one or
two roots in the interval 
$[p\t,1]$, and it has two distinct
roots in this interval 
if and only if $G'(1-p\t)>1$.
\endproclaim

The function $f_p$ is sketched in Figure 
\refer{2.fig1}. Next is the main result
of this section.

\topinsert
\figure
\mletter{y}{2.9}{0.6}
\mletter{\t}{2.9}{4.6}
\mletter{p\t}{4.3}{6.8}
\mletter{y=x}{6}{2.8}
\mletter{y=f_p(x)}{7}{4.4}
\mletter{1}{8.6}{6.8}
\lastletter{x}{9.5}{6.8}
\centerline{\epsfxsize=8cm
  \epsfbox{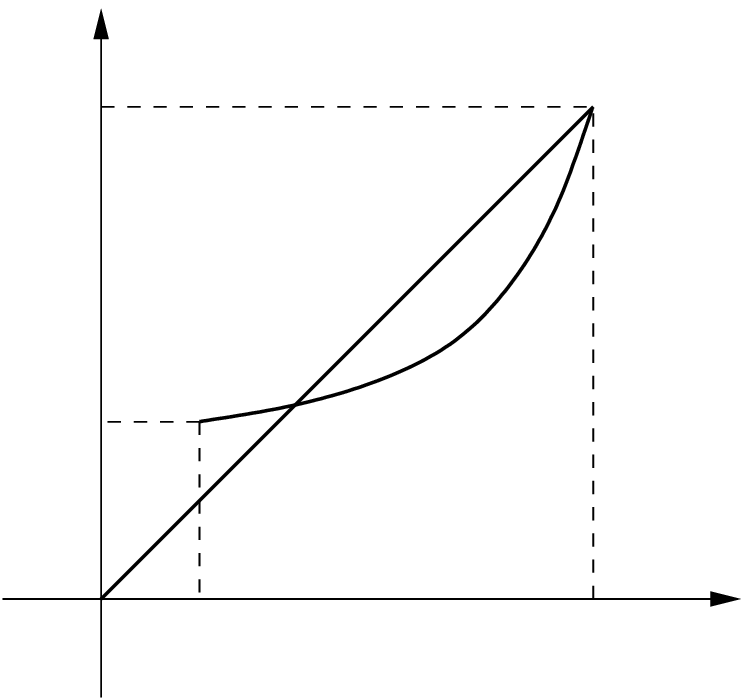}}

\capt{A sketch of the graphs of $y=x$ and 
$y=f_p(x)$
on the interval $[p\t,1]$. Whether or not there exists
a root of $\a=f_p(\a)$ other than at $\a=1$ depends on the
gradient of $f_p$ at $\a=1$.}
\label{2.fig1}
\endinsert

\comment
Let $T=(V,E)$ be given and let $x\in V$. 
The colour of $x$ depends only on the states 
of edges belonging to the $x$-rays of $T$, and thus
the states of edges of $T$ which lie in no ray are immaterial
to whether or not any given vertex is \green. We may
thus delete every edge of $T$ which lies in no ray.  This
amounts to deleting from $T$ all vertices $z$ for which 
the sub-tree rooted at $z$ is finite. 
The effect of this is to replace $T$ by a sub-tree $T'$, and
it well-known (ref?) that $T'$ is itself the family-tree of a branching process
with the property that every individual has a {\it non-empty\/}
family of offspring. We may therefore assume without
loss of generality that $T$ itself has this property, and thus we 
assume henceforth that  the \pgf\ of a typical family-size of $T$
satisfies
$$
G(0)=0, \q 1<G'(1)<\oo.
\TAG
$$\label{2.ass}
\endcomment

\proclaim{Theorem \ch.2}
Consider a branching process whose family-size \pgf\ satisfies \refer{2.ass},
and let $p\in[0,1)$.
Then  $\g$ is the smallest root in the interval
$[p\t,1]$ of the equation
$$
\g = \t+G(\g-p\t).
\TAG
$$\label{2.equ}
The following are equivalent{\rm:}
\ii{\rm(i)} $\g=1$, which is to say there exist almost 
surely no red vertices,
\ii{\rm(ii)} $G'(1-p\t)\le1$, 
\ii{\rm(iii)} $\pga \le p\le 1$, where $\pga$ is given
uniquely by $G'(1-\pga \t(\pga))=1$.
\endproclaim

We make a remark about the value $\pga$ given in 
Theorem 2.2(iii). The function $\t=\t(p)$ is smooth
when $p>p_0=1/G'(1)$. On differentiating \refer{2.extcond2} we find that
$$
\t' = (\t+p\t')G'(1-p\t).
$$
Hence, for $p_0<p < 1$, $G'(1-p\t) \le 1$
if and only if $\t'\le \t+p\t'$, which is to 
say that 
$$
\frac d{dp}\bigl((1-p)\t(p)\bigr) =(1-p)\t' - \t \le 0.
$$
Therefore, $\pga$ is characterised as the value of $p
\in[0,1]$ which maximises $(1-p)\t(p)$. 

We point out that
the coloured tree $T$ constitutes a multi-type branching
process. That is, suppose that each vertex
of $T$  is coloured red, blue, or yellow in
the manner described above. 
We may think of the colour
of any given vertex as its {\it type\/}, 
and then it
is an exercise in the theory of multi-type branching processes to show
that $T$ (when coloured) has the same 
distribution as the 
family-tree of a multi-type 
branching process with
certain offspring-type distributions. This
is a consequence of a general result 
for multi-type processes which may already be known, 
and whose details
are contained in the next section.

We turn briefly to a particular
instance of importance for the \rc\ model on a regular tree.
Let $m\in\{2,3,\dots\}$ be a given integer, and
let $T_m$ denote the infinite labelled 
rooted tree in which the root has degree $m$ and every other vertex has degree $m+1$. This is the family-tree
of the branching process with \pgf\ given by 
$G(\a)=\a^m$, $\a\in\R$. We 
consider bond percolation on $T_m$
with edge-density $p$ as above, and we arrive at the following result
to be found in [\ci{Jon99}].

\proclaim{Corollary \ch.3 [\ci{Jon99}]}
Let $m\in\{2,3,\dots\}$ and $p\in[0,1]$. The probability $\g$ that the 
root of
$T_m$ is \green\ satisfies
$\g=1$ if and only if $p\ge \pg(m)$, where
$$
\pg(m) = \frac{1-m^{-1/(m-1)}}{1-m^{-m/(m-1)}}.
\TAG
$$\label{2.exactpg}
\endproclaim

In particular, 
$$
\pg(2) = \tfrac23,\q \pg(3) = \tfrac3{13}(4-\sqrt 3) 
\approx 0.52337\ldots,
\TAG
$$\label{2.smallm}
and it is easily seen that
$$
\pg(m) \sim \frac{\log m}m \qq\text{as } m\to\oo.
\TAG
$$\label{2.largem}

Finally prior to the proofs, we make a remark about the situation when \refer{2.ass} is not assumed in its entirety, but only 
that $1<G'(1)<\oo$. The branching process is then supercritical, but may be finite
with a strictly positive probability.
Even if $T$ is infinite, it will 
generally contain
vertices $x$ for which the set $\sR_x$ of
$x$-rays is empty, and such vertices are automatically assigned the colour \yellow,
following the rules given towards 
the start of this section. The conclusion
of Theorem 2.2 is easily seen to be
valid in this more general setting.

\demo{Proof of Proposition \ch.1}
The function $f_p$ is non-decreasing and 
strictly convex on the interval $[p\t,1]$.
It is clear by \refer{2.extcond2}
that $f_p(1)=1$ and 
$f_p(p\t)\ge p\t$.
See Figure \refer{2.fig1}.

If $f_p'(1)>1$, there exist two distinct roots of the equation $\a=f_p(\a)$ in
$[p\t,1]$; if $f_p'(1)\le 1$, then $\a=1$ is the unique such root.
\qed\enddemo

\demo{Proof of Theorem \ch.2}
Let $k\ge 1$, and let
$X$ be the number of offspring of the progenitor $0$.
We say that a vertex $x$ is $k$-\yellow\ if: $x$ is not blue,
but every $x$-ray contains a blue vertex belonging to the first $k$
generations of $T$. Vertex $x$ is called 
$k$-\green\ if it
is either blue or $k$-yellow.
Let $\g(k,X)$ be the (conditional) probability given $X$ that
0 is $k$-\green, and write $\g(k)=\PP(\g(k,X))$.
Now, $0$ is $k$-\yellow\ if and only if it is not blue but
every child is either $(k-1)$-\yellow\ or blue. This occurs if and only if
every child $y$ of 0 satisfies:
either $y$ is $(k-1)$-\yellow, or $y$ is blue and the edge $\langle 0,y\rangle$
is closed. Therefore,
$$
[\g(k,X)-\t]  =  
\bigl([\g(k-1)-\t] + (1-p)\t\bigr)^X. 
$$
We take expectations to find that
$$
\g(k)=f_p(\g(k-1))
\TAG
$$\label{2.recur2}
where $f_p$ is given in \refer{2.defnf}. Now $\g(k)\to\g$ as $k\to\oo$,
and $f_p$ is continuous, whence $\g$ satisfies \refer{2.equ}. Since $\g(0)=\t<1$, we have by the usual recursion argument that
$\g$ is the smallest root in $[p\t,1]$ of 
\refer{2.equ}. 
Note that $\g\ge\t$ by \refer{2.defg5}.

By Proposition \ch.1,
$\g=1$ if and only if $f_p'(1)=G'(1-p\t)\le 1$. Now $G$ is strictly convex
and differentiable on $[0,1]$, and $G'(1)>1$ by \refer{2.ass}, while $G'(0)=\PP(X=1)<1$.
Therefore there exists a unique $\b\in(0,1)$ such that $G'(\b)=1$,
and $G'(1-p\t)\le 1$ if and only if 
$1-p\t \le \b$.
\qed\enddemo

\demo{Proof of Corollary \ch.3}
Since $G(\a)=\a^m$, the unique root of the equation $G'(\b)=1$ is given by
$\b=m^{-1/(m-1)}$. If $p\t=1-\b$ then, by
\refer{2.extcond2},
$$
1-\t = G(\b) = \b^m,
$$
whence
$$
p=\frac{1-\b}{\t}=
\frac{1-m^{-1/(m-1)}}{1-m^{-m/(m-1)}}
$$
as claimed.
\qed\enddemo

\endsection

\def\ch{3}
\section{\ch. Multi-type branching processes}
We prove a general result about multi-type processes in this section,
and then apply it 
to the coloured branching processes of Section 2.
A related argument underpins the Markovian construction
of \rc\ measures on trees in [\ci{Hag96}].

Consider a multi-type (Galton--Watson) 
branching process with a 
set $\sI$ of types; $\sI$ may be finite or countably infinite.
We assume, for convenience, that the children of each individual are
ordered in some manner, and we may if
necessary impose a random ordering within
families.
Suppose also that we are given a (measurable) classification of the
possible family-trees into a (finite or countable) set of types $\sJ$
that we call `colours' (not to be confused with the original types that
will be called `types').
We colour each vertex $x$ of the family-tree by the colour of the sub-tree
rooted at $x$.

We shall assume that the colouring rule has the 
property that the colour of any given vertex
is determined by the number and types and colours of its children. If $x$ has $k$ offspring 
labelled $1,2,\dots,k$, the $r$th of
which has type $i_r$ and colour $j_r$, we 
denote the colour of $x$ by
$\Phik$.

\proclaim{Theorem \ch.1}
The coloured family-tree of the process, with vertices marked by both type and
colour, is a multi-type branching process with type space $\sI\times\sJ$.
\endproclaim

\demo{Proof}
Let $\pik$ be the probability that an 
individual of type $i\in\sI$ has
$k\ge0$ children of types $\ik$, respectively.
Let $\qij$ be the probability that the family tree, starting with an
individual of type $i\in\sI$, receives 
colour $j\in\sJ$.

Let
$$
\multline
  \pijk
\\
  =\cases \qij^{-1} \pik\prod_{l=1}^k\qijl 
&\text{if } \Phik=j,\\    
0 &\text{otherwise}.
  \endcases
\endmultline
$$
There is probability 
$\pik\prod_{l=1}^k\qijl$ that a family-tree starting with an individual of
type $i$ has $k$ children with types and colours $\ijk$.
In this case the root is coloured $\Phik$,  
and thus the probability of this
happening and the root being coloured $j$ is $\qij\pijk$.
Consequently, the conditional probability that this happens given that
the root has type $i$ and colour $j$ is $\pijk$.

Moreover, if we label the $k$ children of the root 
as $x_1,\dots,x_k$, and we require
some further events $\sE_1,\dots,\sE_k$, where $\sE_l$ depends only on
the sub-tree rooted in $x_l$, then in
these probabilities we have to replace $\qijl$ by
$$
\multline
\hskip1cm\PP(\sE_l\text{, and $x_l$ receives colour $j_l$}
\mid \text{$x_l$ has type
  $i_l$})
\\  
=\qijl \PP(\sE_l\mid\text{$x_l$ has type $i_l$ and colour $j_l$}).
\hskip1cm
\endmultline
$$
Thus the probabilities 
are multiplied by
$$
\prod_{l=1}^k\PP(\sE_l\mid\text{$x_l$ has type $i_l$ and colour $j_l$}).
$$
Consequently, 
given the types and colours in the first generation,
the $k$ branches are independent of each other and are copies of the
entire coloured family-tree, with appropriate initial conditions.
In other words, the coloured family-tree is a multi-type branching
process with type space $\sI\times\sJ$, where
the probability that a particle of type and colour $(i,j)\in\sI\times\sJ$ has
$k\ge0$ children of types and colours $\ijk$, respectively, is $\pijk$.
\qed\enddemo

Let us apply this result to the coloured trees
of Section 2. The tree $T$ is the family-tree of a branching process with \pgf\ $G$
satisfying \refer{2.ass}. 
We designate any given edge
of $T$ as `open' with probability $p$,
and as `closed' otherwise. Thus the `type' of
any vertex of $T$ is taken to be the state of the incoming edge,
namely either $1$ (open) or $0$ (closed).
(The type of the root of the tree is irrelevant.) The number of offspring of $x$ thus does not depend
on the type of $x$, and each child is assigned a type independently of the types of the other offspring, with probability $\pi_0=1-p$ (\resp, $\pi_1=p$) for type 0 (\resp, type 1).
We
may thus write $q_j$ for the
probability that the root receives colour 
$j\in\set{b,y,r}$.

Consider for example the binary tree $T_2$ defined
at the end of Section 2. We
have that
$$
\align
p_b(i_1,j_1,i_2,j_2)&=
\cases
p\pi_{i_2}q_{j_2}  
&\text{if } (i_1,j_1)=(1,b),
\\
p\pi_{i_1}q_{j_1}  
&\text{if } (i_2,j_2)=(1,b),
\\
0
&\text{otherwise};
\endcases
\\
p_y(i_1,j_1,i_2,j_2)
&=
\cases
\pi_{i_1}\pi_{i_2}q_{j_2}  
&
\text{if $j_1=y$ and either $j_2=y$ or } (i_2,j_2)=(0,b),
\\
\pi_{0}\pi_{i_2}q_{b}  
&
\text{if $(i_1,j_1)=(0,b)$ and } j_2=y,
\\
\pi_{0}^2q_{b}^2/q_{y}  
&
\text{if } (i_1,j_1)=(i_2,j_2)=(0,b),
\\
0
&
\text{otherwise};
\endcases
\\
p_r(i_1,j_1,i_2,j_2)
&=
\cases
\pi_{i_1}\pi_{i_2}q_{j_2}  
&
\text{if $j_1=r$ and } (i_2,j_2)\neq(1,b),
\\
\pi_{i_1}\pi_{i_2}q_{j_1}  
&
\text{if $j_2=r$ and } (i_1,j_1)\neq(1,b),
\\
0
&
\text{otherwise.}
\endcases
\endalign
$$

The following is easily seen.
The mean number of red offspring of a red individual is 
$$
\mu_r=2(1-pq_b)=2(1-p\t),
$$
and $\mu_r>1$
if and only if $p<\pg(2)=\frac23$; see 
\refer{2.exactpg} and Theorem 2.2. 
Therefore, the red sub-tree
of the multi-type process having red progenitor is supercritical if and only if $p<\pg(2)$.
A similar calculation is valid for $T_m$ 
with $m\ge 2$.

In fact, for any branching process as in Section 2, 
a simple
calculation yields that the mean number of red offspring
of a red individual is $G'(1-p\t)$.
\endsection

\def\ch{4}
\section{\ch. Random-cluster measures on trees}
Henceforth, we restrict ourselves to 
the regular infinite labelled 
$m$-ary tree $T_m'=(V,E)$, where $m\in \{2,3,\dots\}$. Each vertex has degree $m+1$,
and there is a 
distinguished {\it origin\/} labelled 0. We
shall state our results for general $m$, but may sometimes consider
the special case $m=2$ for simplicity. Part of $T_2'$ is drawn in Figure
\refer{3.fig1}. 

The tree $T_m'$ differs from $T_m$ only in
the degree of its root. 
We have chosen to work with $T_m'$
rather than $T_m$ in this
section only because this
is the more natural setting for the 
\rc\ model. Since $T_m'$ is a regular tree,
it has a larger family of graph-automorphisms.
When it comes to calculations of
critical values and the like, the differences between $T_m'$ and $T_m$
are largely cosmetic.

\topinsert
\figure
\lastletter{0}{6}{2.5}
\centerline{\epsfxsize=5cm
  \epsfbox{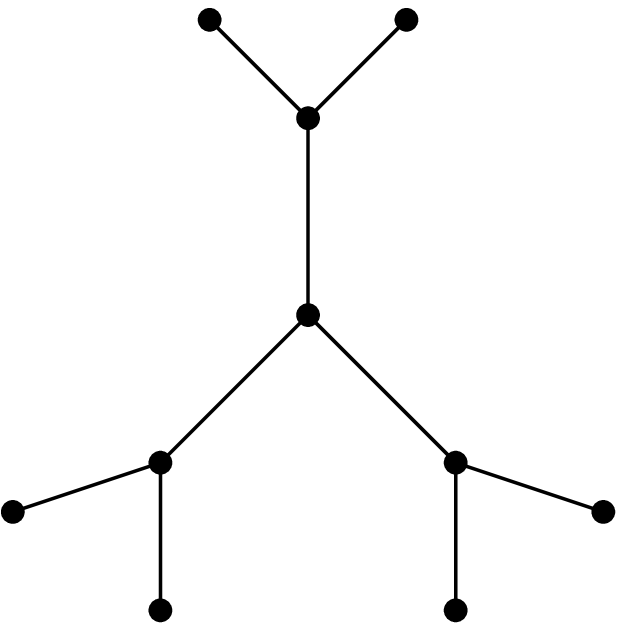}}

\capt{Part of the infinite binary tree $T_2'$.}
\label{3.fig1}
\endinsert

We continue the study of \rc\ measures on $T_m'$ 
initiated in [\ci{Hag96}], beginning with a more general definition
than that used there. Let $\Om=\{0,1\}^E$, and 
equip $\Om$ with the $\s$-field $\sF$ generated by
the finite-dimensional cylinder sets. An edge $e$ is called {\it open\/}
in a configuration $\om$ ($\in\Om$) if $\om(e)=1$,  and
{\it closed\/} otherwise. We write $\eta(\om) = \{e: \om(e)=1\}$
for the set of open edges in $\om$. We shall consider probability
measures on the measurable pair $(\Om,\sF)$ which satisfy a certain
`\rc' condition. 
Since $T_m'$ contains no circuits, \rc\
measures on $T_m'$ 
are simply product measures. A much more
interesting structure is revealed 
through the introduction of the
concept
of `boundary conditions'.  A similar
development for Ising models has been explored in
the statistical physics literature, see [\ci{BRZ}, \ci{BRSSZ}], and in
the probability literature under the title `broadcasting in trees',
see [\ci{Mar}, \ci{Mos03}].
Boundary conditions may be introduced
in the more general context of non-amenable graphs, but we do
not follow this route here; see [\ci{G04}, \ci{HJL}, \ci{Jon99}] for
accounts of the \rc\ model on a non-amenable graph.

Each edge of $T_m'$ is directed away from the root 0. We shall make use of the rays of $T_m'$ and we remind
the reader of the notation concerning 
rays at the beginning of Section 2. The set $\sR$ of 
0-rays
of $T_m'$ is in one--one correspondence with the set 
$\{1,2,\dots,m+1\}\times \{1,2,\dots,m\}^\N$.
Let $\sim$ be an equivalence relation on $\sR$, 
and write $\sC(\sim)$ 
for the family of equivalence classes of $\sim$. Let $\sE$ denote
the class of all equivalence relations on $\sR$. There is a partial order
$\le$ on $\sE$ given by: ${\sim_1} \le {\sim_2}$ if: 
$$
\text{for all $\rho,\rho'\in\sR$, $\rho\sim_2 \rho'$ 
whenever $\rho\sim_1\rho'$.}
\TAG
$$\label{4.partord}
There is a minimal (\resp, maximal) 
equivalence relation in this 
partial order which we denote as $\sim^0$ 
(\resp, $\sim^1$). The equivalence classes of $\sim^0$ are singletons,
whereas $\sC(\sim^1)=\{\sR\}$.

For $x\in V$, let $\Pi_x$  be the set
of infinite (undirected) paths of $T_m'$ with endpoint $x$.
Let  $x\in V$ and let
$\pi\in\Pi_x$; 
there exists a unique 0-ray, denoted
$\rho_\pi$, such that $\pi$ and
$\rho_\pi$ differ on only
finitely many edges. For given $x$, this
gives a one--one correspondence  $\pi\leftrightarrow\rho_\pi$
between $\Pi_x$ and $\sR$. (For $x=0$, it is the identity.)
Any relation
${\sim}$ on $\sR$ may be extended to a
relation on $\bigcup_{v\in V}\Pi_v$ by:
for $\pi\in\xsR_u$, $\tau\in\xsR_v$, we have 
$\pi \sim \tau$ if and only if 
$\rho_\pi\sim\rho_\tau$.

For any vertex $u\in V$, we write $\sR_u'$ for the subset
of $\sR$ comprising all rays which pass through $u$. 
The correspondence between $\Pi_u$ and $\sR$ restricts to a correspondence
between their subsets $\sR_u$ and $\sR_u'$, such that
any ray $\rho_u\in\sR_u$ corresponds to the
unique ray $\rho_u'\in\sR$
of which it is a sub-ray.

The equivalence relation $\sim$ ($\in\sE$)
may serve as a boundary condition,
to be interpreted roughly as follows. Suppose that $\om\in\Om$, and let
$C_1$, $C_2$ be two distinct
connected components of the graph $(V,\eta(\om))$.
Then $C_1$ and $C_2$ are considered to be the same cluster if there
exist rays $\rho_1\subseteq C_1$, $\rho_2\subseteq C_2$ such 
that $\rho_1\sim \rho_2$. Otherwise expressed, two rays are `identified
at infinity' if they are equivalent under $\sim$. This will be made more
rigorous in the following formal definition of a \rc\ measure.

Let $\La$ be a finite subset of $V$, and let $E_\La$ be the set of
edges of $T_m'$ having both endvertices in $\La$. 
For $\xi\in\Om$,
we write $\Om_\Lambda^\xi$ for
the (finite) subset of $\Om$ containing all configurations $\om$ satisfying
$\om (e)=\xi (e)$ for $e\in E\sm E_\La$;
these are the configurations which `agree with $\xi$ off $\La$'.

Let ${\sim}\in\sE$, $\xi\in\Om$, and $\om \in\Om_\La^\xi$. 
The configuration
$\om$ gives rise to an `open graph' on $\La$, namely
$G(\La,\om)=(\La,\eta(\om)\cap E_\La)$. We augment this graph by adding
certain new edges. Specifically, for distinct
$u,v\in\La$,
we add a new edge between the pair $u$, $v$ if either:
\ii{(a)} there exists a path of
$E\sm E_\La$ from $u$ to $v$ which is open in
$\xi$, or
\ii{(b)} there exists
a vertex-disjoint pair of infinite paths $\xrho_u\in \xsR_u$, 
$\xrho_v\in\xsR_v$ 
satisfying $\xrho_u\sim\xrho_v$,
which are open in  $\xi$ and which are 
edge-disjoint from $E_\La$.
\flushpar
We write $G^{\xi, \sim}(\La,\om)$ for the resulting
augmented
graph, and we let $k^{\xi,\sim}(\La,\om)$ be the number of connected components
of $G^{\xi,\sim}(\La,\om)$. These definitions are motivated
by the idea that each equivalence class of rays 
leads to a common `point at infinity',
through which vertices $u$ and $v$ may be 
connected by open paths.

Next we define a \rc\ measure corresponding to
a given equivalence relation $\sim$. 
Let ${\sim}\in\sE$,  $\xi\in\Om$, and let
$p\in[0,1]$ and $q\in(0,\oo)$. 
We define $\phi_{\Lpq}^{\xi,\sim}$
to be the random-cluster measure on the finite graph
$(\Lambda ,E_\Lambda )$ `with boundary condition $\xi$'.
More precisely, let $\phi_{\Lpq}^{\xi,\sim}$
be the probability measure on the pair $(\Om ,\sF )$ given by
$$
\phi_{\Lpq}^{\xi,\sim} (\om ) = \cases\dfrac 1{\Zlpqxs}
\bigg\{\displaystyle\prod_{e\in E_\Lambda} p^{\om (e)} (1-p)^{1-\om (e)}\bigg\}
q^{k^{\xi,\sim}(\La,\om)} &\text{if $\om\in\Om_\Lambda^\xi$},\\
0 &\text{otherwise},
\endcases
\TAG
$$\label{3.defrc}
where
$\Zlpqxs$ is the appropriate normalising constant
$$
\Zlpqxs =\sum_{\om\in\Om_\Lambda^\xi}\bigg\{\prod_{e\in E_\Lambda}
p^{\om (e)} (1-p)^{1-\om (e)}\bigg\} q^{k^{\xi,\sim}(\La,\om)}.
\TAG
$$\label{old3.12}
Note that $\phi_{\Lpq}^{\xi,\sim} (\Om _\La^\xi) = 1$.

In the special case when $\xi=1$ and 
${\sim}={\sim^1}$, we write
$\phi_{\Lpq}^{1}$ for $\phi_{\Lpq}^{\xi,\sim}$. This is usually called
the \rc\ measure on $\La$ with `wired' boundary conditions, and
it has been studied in a slightly disguised form in [\ci{Hag96}].  As
usual (see [\ci{G02}] and the references therein),  
\rc\ measures satisfy the FKG inequality when $q \ge 1$,
and this allows the deduction that the weak limit
$$
\phi_{p,q}^1 = \lim_{\La\uparrow V} \phi_{\Lpq}^{1}
\TAG
$$\label{3.weakl}
exists when $q \ge 1$, and is independent 
of the manner in which
the limit $\La\uparrow V$ is taken. As a side-comment,
we remark that the FKG inequality is a fundamental
technique in the study of the \rc\ model. 
This is already very familiar in the field 
(see [\ci{G02}], for example), and we do
not explain it further here. Thus, we shall omit
details of arguments involving the FKG inequality
and the stochastic ordering relation $\lest$
of probability measures. We note for later use that,
for all $q\ge 1$ and ${\sim}\in \sE$,
$$
\phi_{\De,p,q}^{1,\sim} \lest \phi_{\La,p,q}^{1,\sim}
\qq\text{ if } \De \supseteq \La.
\TAG
$$\label{4.stord} 

For any finite subset  $\La\subseteq V$, 
we write $\sT_\La$ for the $\sigma$-field generated
by the set $\{\om(e):e\in E\sm E_\La\}$ of states of edges
having at least one endvertex outside $\Lambda$.

Let ${\sim}\in\sE$, $0\le p\le 1$, and $q>0$. 
A probability measure $\phi$ on $(\Om ,\sF)$ 
is called a {\it $(\sim)$random-cluster measure\/} with parameters $p$ and
$q$ if: for
all $A\in\sF$ and all finite subsets $\Lambda\subseteq V$,
$$
\phi (A\mid\sT_\Lambda )(\xi)=
\phi_{\Lpq}^{\xi,\sim} (A)\q\text{for $\phi$-a.e.
$\xi$.}
\TAG
$$\label{3.rcdef}
The set of such measures is denoted $\sR^\sim_{p,q}$. Random-cluster
measures were introduced 
in [\ci{G93}] and are studied systematically in 
[\ci{G04}] and elsewhere.
Note that the cases when $p=0$, $p=1$ or $q=1$ are trivial; in these
cases, \refer{3.defrc} gives the product measure on $(\La,E_\La)$ for
every $\xi$ and $\sim$, and then \refer{3.rcdef} holds if and only if
$\phi$ is the product measure $\phi_p$ with density $p$ 
on $(\Om,\sF)$.

This is an appropriate moment to introduce a measurability assumption 
on $\sim$.
The left hand side of \refer{3.rcdef} is a measurable function of
$\xi$, so we want the right hand side to be a measurable function of
$\xi$ too, with respect to $\sF$ or at least with respect to its completion
for $\phi$.
For $\La$ a finite subset of $V$ and $u,v\in\La$,
let $K_{u,v,\La}^\sim$
denote the set of
$\om\in\Om=\{0,1\}^E$ such that there exist
infinite vertex-disjoint
paths $\pi_u\in\Pi_u$, $\pi_v\in\Pi_v$
satisfying $\pi_u\sim\pi_v$ and that are open
in $\om$ and edge-disjoint from $E_\La$.
We call the equivalence relation
$\sim$ {\it measurable\/} if 
$K_{u,v,\La}^\sim\in\sF$ for all such $u,v,\La$.
It is an easy 
exercise to deduce, if $\sim$ is measurable, that
$\phi_{\La,p,q}^{\xi,\sim}(A)$ is
a measurable function of $\xi$, thus permitting
condition \refer{3.rcdef}.
We write $\sEm$ for 
the set of all measurable members of $\sE$.
We discuss measurability further in the next section.

Returning to the extremal equivalence relations $\sim^0$, $\sim^1$,
for simplicity of notation 
we write $\sR^{\sim^0}_{p,q}=\sR^{0}_{p,q}$
and similarly $\sR^{\sim^1}_{p,q}=\sR^{1}_{p,q}$.
Members of $\sR^0_{p,q}$ (\resp, $\sR^1_{p,q}$) are
called `free' \rc\ measures (\resp, `wired' \rc\ measures).  

The basic questions of interest include the following. 
For what ${\sim}\in\sEm$, $p$, $q$ is the set 
$\sR^\sim_{p,q}$  non-empty, and when does 
$\sR^\sim_{p,q}$ comprise a singleton only?
Progress towards answers has been made in [\ci{Hag96},
\ci{Jon99}].
Let $\phi_p$ denote product measure with density $p$ on $(\Om,\sF)$. 
We define $\pi:[0,1]\times(0,\oo)\to[0,1]$ by
$$
\pi(p,q)=\frac p{p+q(1-p)}.
\TAG
$$\label{3.pprime}
When $p$, $q$ are given, we use the abbrevation
$\pi=\pi(p,q)$.
Note that $\pi\ne p$ except in the trivial cases $p=0$, $p=1$ or $q=1$.

\goodbreak
\proclaim{Theorem \ch.1\ [\ci{Hag96}]}
\flushpar
\ii{\rm(a)} For $0\le p\le 1$ and $q>0$, the set $\sR^0_{p,q}$ of free
\rc\ measures 
comprises the singleton $\phi_{\pi}$ only, where 
$\pi=\pi(p,q)$.
\ii{\rm(b)}  For $0\le p\le 1$ and $q> 0$, the set
$\sR^1_{p,q}$ of wired \rc\ measures
is non-empty. If $q\ge 1$, then $\sR^1_{p,q}$ 
contains the weak limit
$\phi_{p,q}^{1}$ given in \refer{3.weakl}.
\endproclaim

We present in Theorem 6.1(a) a necessary and sufficient condition
for the statement $\phi_\pi\in \sR_{p,q}^\sim$, 
for $\sim$ belonging to a certain class of equivalence
relations to be defined in the next section.
 
We write $\{x\lra y\}$ for the set of all $\omega\in\Omega$ for
which there exists an open path joining vertex $x$ to vertex $y$. (For
the moment, we refer to open paths in the {\it undirected\/} 
tree $T_m'$.) The complement
of the event $\{x\lra y\}$ is denoted $\{x\nlra y\}$. We write
$\{x\lra\oo\}$ for the event that $x$ is the endvertex of
an infinite open path of $T_m'$. For  
a probability measure $\phi$ on $(\Om,\sF)$, we define
the {\it percolation probability\/} by
$$
\t(\phi) =\phi(0\lra \oo).
\TAG
$$\label{3.percprob}
Of particular interest are the two special cases
$$
\t^0(p,q) = \t(\phi_{\pi}),\q \t^1(p,q) = \t(\phi_{p,q}^1),
$$
where the former is well-defined for all $p$, $q$, and the latter 
when $q\ge 1$ at least. 

The function $\t^0(p,q)$ is the survival probability of a branching
process with $\bin(m,\pi)$ family-sizes (subject to the small change of
vertex degree at the origin). It may be computed as in
\refer{2.extcond}--\refer{2.extcond2}. In particular, its {\it critical
value\/}
$$
\pc^0(q) = \sup\{p: \t^0(p,q)=0\}
\TAG
$$\label{4.crit1}
is immediately seen to be the value of $p$
for which $\pi(p,q)=m^{-1}$, whence
$$
\pc^0(q) = \frac q{m+q-1}.
\TAG
$$\label{4.crit2}
Less standard is the calculation given in 
[\ci{Hag96}] 
of the critical value
$$
\pc^1(q) = \sup\{p: \t^1(p,q)=0\}
$$
when $q \ge 1$, namely
$$
\pc^1(q) = \cases \pc^0(q) &\text{if } 1\le q\le 2,\\
U_q &\text{if $q > 2$},
\endcases
$$
where $U_q$ is the unique value of $p \in (0,1)$ for which the
polynomial
$$
(q-1)x^{m+1} + \left(1-\frac p{1-p} - q\right)x^m +\frac 1{1-p} x - 1
$$
has a double root in $(0,1)$.  Applying this as in [\ci{Hag96}]
when $m=2$, we find that
$$
\pc^1(q)=\cases \dfrac q{q+1} &\text{if } 1\le q\le 2,\\
\dfrac{2\sqrt{q-1}}{1+2\sqrt{q-1}} &\text{if } q > 2.
\endcases
\TAG
$$\label{4.pc1}
[An alternative proof of these facts
may be obtained by the parallel/series 
replacement method used in the proof
of Theorem 6.4.]

We note from [\ci{Hag96}] for later use that,
for $q\ge 1$, 
$$
\t^1(p,q) > 0 \q\text{if and only if}\q
\cases 1\le q\le 2,\ p>\pc^1(q),\ \text{or}\\
q > 2,\ p \ge \pc^1(q).
\endcases
\TAG
$$\label{4.I}

Let $q\ge 1$.
It was proved in [\ci{Hag96}] that there exists a continuum of
probability measures in $\sR^1_{p,q}$ when $q> 2$
and $\pc^1(q)\le p\le \pc^0(q)$,
and it was conjectured that $\sR^1_{p,q}$ contains exactly
one measure when $p>\pc^0(q)$. Uniqueness was
proved in [\ci{Jon99}] for sufficiently large values of $p$,
and we recall this result next.
In the notation of Section 2, we take as `mother
process' the process in which every individual has exactly $m$ children,
the corresponding family-size \pgf\ is given by $G(\a)=\a^m$, $\a \in \R$.
On this graph, we construct bond percolation with density $\pi$ given
in \refer{3.pprime}, and we
ask for the probability that the origin is \green.
By Corollary 2.3 (see also [\ci{Jon99}]), 
the probability $\g$ that the origin is \green\  
satisfies
$$
\g=1\q\text{if and only if} \q \pi \ge \pg(m),
$$
where $\pg(m)$ is given in \refer{2.exactpg}.

\proclaim{Theorem \ch.2 [\ci{Jon99}]} Let $q\ge 1$ and let $p$
be such that 
$$
\pi=\frac p{p+q(1-p)}\ge \pg(m).
$$
The set
$\sR^1_{p,q}$ comprises the singleton $\phi^1_{p,q}$ only.
\endproclaim

\demo{Proof}
This is a special case of the forthcoming Theorem 6.3.
\qed\enddemo

\endsection

\def\ch{5}
\section{\ch. Relations on the set of rays}
We consider next the case of a general boundary condition $\sim$ ($\in\sEm$).  
We cannot prove in general that 
$\sR^{\sim}_{p,q}$ is non-empty, 
but only for a certain class of equivalence relations which
we introduce in this section. 
It is in fact unnecessary for the
present purposes to require a relation
on $\sR$ to be an {\it equivalence\/} relation, and thus
we shall broaden the discussion in this section
to the class of all 
symmetric relations on $\sR$. 

For simplicity, we continue to consider 
the infinite $m$-ary tree $T_m'$ with root 0,
and every edge oriented away from 0. The conclusions of
this section are valid under considerably less restrictive
assumptions on the underlying tree.

We write $\sSR$ for
the class of all symmetric reflexive relations on $\sR$, and
we think of a relation as a subset of the space $\sR^2$.
Thus we consider the set $\sSR$ of subsets $S$ of $\sR^2$ such that:
\ii{(a)} $(\rho,\rho)\in S$ for all $\rho\in \sR$,
\ii{(b)} $(\rho_1,\rho_2)\in S$ whenever $(\rho_2,\rho_1)\in S$.

The set $\sR$ of rays may be viewed
as a topological space
with the product topology.  Since $\sR$ is
the product of compact spaces, it is itself compact. 
The  family $\{\sR_v': v\in V\}$ forms a base
for the space. 
We endow $\sR^2$ with the product topology, and we call 
the relation
$S\in\sSR$ {\it closed\/} (\resp, {\it open, Borel\/})
if $S$ is a closed (\resp, open, Borel) subset of 
$\sR^2$. We write $\sSRc$ (\resp, $\sSRf$, $\sSRm$)
for the set of closed (\resp, open, Borel)
relations. By definition, $\sSRc,\sSRf\subseteq \sSRm$.
We shall see in Corollary 5.6 that every open 
equivalence relation is closed.

There follows a description of a certain family
of closed relations.
A (finite or infinite)
subset $W$ of $V$ is called {\it incomparable\/} if there exists
no 0-ray of $T_m'$ which contains more than one member of $W$; $W$ is
called a {\it cutset\/} if it is incomparable 
and every 0-ray contains some
member of $W$. 
The smallest cutset is the singleton
set $\{0\}$. Let $W$ be a cutset.  Each $\sR_w'$, $w\in W$,
is an open subset of the compact space $\sR$, and in
addition $(\sR_w': w\in W)$
is a cover of $\sR$. By compactness, $W$ is finite. In summary,
all cutsets are finite. Every cutset $W$ generates 
a finite set of vertices $\ol W = W \cup \ins(W)$, where
$\ins(W)$ is 
the set of all vertices
$v\in V$ such that there exists a path of $T$ from $0$ to $v$ which is
vertex-disjoint 
 from $W$. We refer to such 
a set $\ol W$ as a
{\it box\/} of the tree.

Let $W$ be an incomparable (finite or infinite) set
of vertices, and partition $W$ as $W=W_1\cup W_2$.
Let $S$ be the (equivalence) relation given by
$$
S = \biggl(\bigcup_{w\in W_1} \sR_w'\times \sR_w'\biggr)
\cup  \biggl(\bigcup_{w\in W_2} \{(\rho,\rho): \rho\in \sR'_w\} \biggr)
\cup \Bigl(\crw\times\crw\Bigr)
$$
where
$$
\crw=\sR\sm \biggl(\bigcup_{w\in W} \sR_w'\biggr).
$$
In the usual jargon borrowed from 
the theory of electrical networks, the relation $S$
corresponds to `wired boundary conditions'
on $\sR_w'$ for every $w\in W_1$, `free boundary conditions'
on $\sR_w'$ for every $w\in W_2$, and a `wired boundary condition'
on the union of all other rays. Certainly $S$ is
an equivalence relation, and in addition
$S$ is a closed relation (it may be considered easier to see
that $\sR^2 \sm S$ is open). The construction of 
[\ci{Hag96}, Section 5] gives rise to equivalence
relations of the above type.
Note that the minimal and maximal relations $\sim^0$ and $\sim^1$
are both of this type (for
$\sim^1$, take $W=\{0\}$); in particular, they are closed.

A simple example of the above recipe
arises when $W$ is a cutset, and is therefore finite.
Let $S^W$ be the relation: 
$$
S^W = \bigcup_{w\in W} \sR_w'\times \sR_w'.
$$
The
equivalence classes of $S^W$
are the sets $\sR_w'$, $w\in W$.
We call $S$ ($\in\sSR$) a {\it cutset relation\/}
if there exists a 
cutset $W$ such that $S=S^W$. 
The maximal
equivalence relation $\sim^1$ is a cutset relation with
single equivalence class $\sR$.

We next continue the measurability discussions from Section 4.
Recall first that, if
$x\in V$ and 
$\pi$ is an infinite path of $T$ with endvertex $x$, then 
there exists a unique 0-ray, denoted 
$\rho_\pi$, such that $\pi$ and 
$\rho_\pi$ differ on only
finitely many edges. 
(The path $\pi\in\Pi_x$ comprises
a finite path from $x$ to some vertex $w$ directed
in the opposite direction, followed by a ray in 
$\sR_w$ which differs from $\rho_{\pi}$ only in the 
absence of the first section of $\rho_\pi$ from 0 to $w$.)
For given $x$, the pair $\pi$, $\rho_\pi$
are in one--one correspondence. For $\pi$ a ray in $\sR_x$, we 
denote  $\rho_\pi$ by $\pi'$.
We consider five related definitions of subsets of the configuration
space $\Om$.
Let $\La$ be a finite subset of $V$.
\define\asx#1#2{A^S_{#1,#2}}
\define\asuv{\asx uv}
\define\asuvi{A^{S'}_{u,v}}
\define\pp#1#2{\pi_{#1,#2}}
\define\ppp#1#2{P_{#1,#2}}
\ii{1.} $\asuv$, for $u,v\in V$, is the set of configurations such
that there exist open rays $\rho_u\in\sR_u$, $\rho_v\in\sR_v$
with $(\rho'_u,\rho'_v)\in S$.
(These rays have endpoints $u$, $v$.)
\ii{2.} $\asuvi$, for $u,v\in V$, is the set of configurations such
that there exist open rays $\rho_u\in\sR'_u$, $\rho_v\in\sR'_v$
with $(\rho_u,\rho_v)\in S$.
(These rays have common endpoint $0$.)
\ii{3.} $K^S_{u,v,\La}$, for 
$u,v\in\La$, is 
 the subset of $\Om$ such that there
exist infinite open
vertex-disjoint paths $\pi(u)\in\Pi_u$, $\pi(v)\in\Pi_v$ satisfying
$(\rho_{\pi(u)},\rho_{\pi(v)})\in S$ and that are
edge-disjoint from $E_\La$.
\ii{4.} $K^S_{u,v}$, for $u,v\in V$, is similarly defined as the subset
of $\Om$ such that there exist vertex-disjoint open paths
$\xrho(u)\in\Pi_u$, $\xrho(v)\in\Pi_v$
satsfying $(\rho_{\xrho(u)},\rho_{\xrho(v)})\in S$. Thus, 
$K_{u,v}^S = K_{u,v,\{u,v\}}^S$.  
\ii{5.} $K^S_{e}$, for an edge 
$e=\langle u,v\rangle \in E$, 
equals $K^S_{u,v}=K^S_{u,v,\{u,v\}}$.

We say that two vertices $u,v$ are {\it comparable\/} if one
is on the path from 0 to the other (including $u=0$, $v=0$ and $u=v$),
and {\it incomparable\/} otherwise.

\proclaim{Theorem \ch.1} 
The following are equivalent, for every $S\in\sSR$.
\ii{\rm(i)} $\asuv\in\sF$ for all $u,v\in V$.
\ii{\rm(ii)} $\asuvi\in\sF$ for all $u,v\in V$.
\ii{\rm(iii)} $\asuv\in\sF$ for all incomparable $u,v\in V$.
\ii{\rm(iv)} $K^S_{u,v}\in\sF$ for all $u,v\in V$.
\ii{\rm(v)} $K^S_{u,v,\La}\in\sF$ for all finite subsets $\La\subseteq V$ and
$u,v\in\La$.
\ii{\rm(vi)} $K^S_{e}\in\sF$ for all edges 
$e\in E$.
\endproclaim

\demo{Proof}
For two vertices $x,y\in V$, let $\pp xy$ be the path between $x$ and
$y$, and let $\ppp xy$ be the cylinder event that all edges in
$\pp xy$ are open.

(i)$\implies$(ii):
$\asuvi=\asuv\cap\ppp0u\cap\ppp0v$.

(ii)$\implies$(iii):
For $\om\in\Om$, let $\phi(\om)$ equal $\om$ except that all edges in 
$\pp0u$ and $\pp0v$ are open. Then $\phi:\Om\to\Om$ is measurable, and
$\asuv=\phi^{-1}(\asuvi)$.

(iii)$\implies$(i):
Suppose that $u$ lies between $0$ and $v$. Let $W$ be the finite set
of the possible vertices $w$ where a ray $\rho_u$ from $u$ that does
not pass $v$ may first
leave the path $\pp uv$, 
Then
$\asuv$ is the finite union of
$\ppp uw\cap\asx wv$, $w\in W$, and $\ppp uv\cap B_v$,
where $B_v$ is the event that
there is an infinite open ray in $\sR_v$.

(i)$\implies$(iv):
Considering the places where the paths $\pi_u$ and $\pi_v$ leave 
$\pp0u$ and $\pp0v$, we see that $K^S_{u,v}$ equals a finite union of
sets $\ppp ux \cap \ppp vy \cap \asx xy$.

(iv)$\implies$(v): Define $\phi(\om)$ to equal $\om$ except that all
edges in $E_\La$ are deemed closed. Then 
$K^S_{u,v,\La} = \phi^{-1}(K^S_{u,v})$.

(v)$\implies$(vi):
A special case.

(vi)$\implies$(iii):
Let $w$ be the father of $v$. 
Define $\phi(\om)$ to equal $\om$ except that all edges in $\pp uw$
are deemed open, and all edges that are incident
 to this path, except at $u$,
are closed. Then
$\asuv = \phi^{-1}(K^S_{\la w,v\ra})$.
\qed\enddemo

We say that the relation $S$ is {\it measurable\/} if the equivalent
conditions in Theorem 5.1 are satisfied.
Note that this agrees with the previous use of 
the the word `measurable' as applied to 
equivalence relations.
We now investigate measurability further. Let 
$O\subseteq\sR\times\Om$
 be the set of pairs $(\rho,\om)$ such that
$\rho$ is a 0-ray that is open in $\om$; this is a closed subset of
the compact space $\sR\times\Om$. 
Further, let
$D\subseteq\sR\times\sR\times\Om$ be the set of triples
$(\rho_1,\rho_2,\om)$ such that both $\rho_1$ and $\rho_2$ are open in
$\om$, and $(\rho_1,\rho_2)\in S$. Thus,
$$
D=(\sR\times O)\cap(\widetilde{\sR\times O})\cap(S\times\Om),
$$
where $\widetilde A=\{(\rho_1,\rho_2,\om):(\rho_2,\rho_1,\om)\in A\}$.
Hence, $D$ is a closed, and thus compact, subset of
$\sR\times\sR\times\Om$ if $S$ is closed, and $D$ is Borel if $S$ is
Borel.
We can now state the following links between the
measurability of a  relation $S$ and the properties of $S$ 
viewed as a subset of $\sR^2$.

\proclaim{Theorem \ch.2} 
\ii{\rm(a)} A closed relation is measurable.
\ii{\rm(b)} A measurable relation is Borel.
\endproclaim

\demo{Proof}
(a)
If $u,v\in V$, then
$\asuvi=\pi_3\bigl(D\cap(\sR'_u\times\sR'_v\times\Om)\bigr)$, where
$\pi_3$ denotes the projection on the third factor. If $S$ is closed,
this is the projection of a compact set, and thus compact, so (ii) in
Theorem 5.1 is satisfied.
(This part is also an immediate consequence of Theorem \ch.4 below.)

(b) Let $\psi:\sR^2\to\Om$ be defined by $\psi(\rho_1,\rho_2)(e)=1$
if and only if $e\in\rho_1\cup\rho_2$. In other words,
$\psi(\rho_1,\rho_2)$ is the configuration with all edges in $\rho_1$
and $\rho_2$ open, but no others. The function
$\psi$ is continuous, and thus
(Borel) measurable.

Let $S$ be a measurable relation, so that each of the six
parts of Theorem 5.1 is valid.
Two distinct 0-rays $\rho_1$ and $\rho_2$ pass though two incomparable
vertices $u$ and $v$, and they satisfy $(\rho_1,\rho_2)\in S$
if and only if $\rho_1\in\sR'_u$, $\rho_2\in\sR'_v$, and
$\psi(\rho_1,\rho_2)\in\asuvi$. Hence,
$$
S=\biggl\{\bigcup_{u,v\,\text{incomparable}} (\sR'_u\times\sR'_v)\cap\psi^{-1}(\asuvi)
\biggr\}
\cup\bigl\{(\rho,\rho):\rho\in\sR\bigr\},
$$
which is a Borel subset of $\sR^2$, using (ii) in Theorem 5.1.
\qed\enddemo

Theorem 5.2 leaves an obvious gap.
\proclaim{Problem \ch.3} 
Is every Borel relation measurable?
\endproclaim

Note that the proof of Theorem 5.2(a) breaks down for Borel relations
because the projection of a Borel set is in general not a Borel set.
However, for Polish spaces (and the spaces we consider are such), the
projection of a Borel set 
is a Suslin (or analytic) set, and such sets are universally
measurable, that is,
they belong to the completion $\sF_\mu$
of $\sF$ for every finite Borel measure $\mu$.
(See, for example, [\ci{Cohn}, Sections 8.2, 8.4].)
Theorem 5.1 holds also if we replace $\sF$ by the $\sigma$-field $\bigcap_\mu
\sF_\mu$ of universally measurable sets, and we thus see that every
Borel relation is (at least) universally measurable, which is enough for
the definition \refer{3.rcdef} to make sense.

The closed relations have a certain property
that will enable a large-volume limit for \rc\ measures,
and we present this next. 
For any box $\La$ and
any $\om\in\Om$, we write 
$\om_\La^1$ for the configuration
that agrees with $\om$ on $E_\La$ and equals 1 
elsewhere, which
is to say that
$$
\om_\La^1(e)=\cases \om(e) &\text{if } e\in E_\La,\\
1 &\text{otherwise}.
\endcases
$$
We define the event 
$$
\K_{e,\La}^S=\{\om\in\Om: \om_\La^1\in K_e^S\},
$$
i.e., if $e=\langle u,v\rangle$, the set 
of configurations such that there exist vertex-disjoint paths
$\xrho(u)\in\Pi_u$, $\xrho(v)\in\Pi_v$
with $(\rho_{\xrho(u)},\rho_{\xrho(v)})\in S$, such that the parts of
the paths inside $\La$ are open.
Note that $\K_{e,\La}^S$ is a cylinder event, and that
it is decreasing in $\La$.

\proclaim{Theorem \ch.4} Let $S \in \sSR$. We have that
$$
\text{for all $e\in E$,}\qq \K_{e,\La}^S \downarrow
 K_e^S\qq\text{as } \La\uparrow V,
\TAG
$$\label{6.simcl}
if and only if $S$ is closed.
\endproclaim

\demo{Proof} 
Assume first that $S$ is closed.
Let $e=\langle x,y\rangle \in E$ 
and $\om\in \Om$.

Let $\Pi_x^e$ (\resp, $\Pi_y^e$) be the set
of infinite (undirected)
paths of $E\sm\{e\}$ with endpoint $x$ (\resp, $y$), and
let $F_\La=F_\La(\om)$ be the set of
all pairs $(\rho_1,\rho_2)\in S$ such that:
\ii{(i)} $\rho_1=\rho_{\nu(x)}$ for some $\nu(x)\in
\Pi_x^e$,
\ii{(ii)} $\rho_2=\rho_{\nu(y)}$ for some $\nu(y)\in
\Pi_y^e$,
\ii{(iii)} all edges in $\nu(x)$, $\nu(y)$ which belong to 
$E_\La$ are open.
\flushpar
Then $F_\La$ is the intersection of $S$ with
a product of two closed sets of rays, and is therefore 
closed and hence compact.
Furthermore, $F_\La$ is decreasing in $\La$. 

We similarly define $F=F(\om)$ by (i) and (ii), but replacing (iii) by:
\ii{(iii$'$)} all edges in $\nu(x)$, $\nu(y)$ are open.

Since $\nu(x)$ and $\nu(y)$ are
uniquely determined by $\rho_1$ and $\rho_2$, it is clear that
$$
F= \bigcap_{\La} F_\La.
$$
Since the sets are compact, this implies that
$$
F(\om)\neq\es \iff F_\La(\om)\neq\es \text{ for every $\La$}.
\TAG
$$\label{sofie}
We now observe that 
$\om\in \K^S_{e,\La} \iff F_\La(\om)\neq\es$, and
$\om\in K^S_{e} \iff F(\om)\neq\es$.
Thus \refer{sofie} can be written 
$K^S_e = \bigcap_\La \K^S_{e,\La}$, which is 
\refer{6.simcl}.

Suppose conversely that \refer{6.simcl} holds.
Let $(\rho_1^n,\rho_2^n)$, $n\ge 1$, be a sequence in $S$ such that
$$
(\rho_1^n,\rho_2^n) \to (\rho_1,\rho_2)\qq
\text{as } n\to\oo
\TAG
$$\label{6.rholim}
for some $\rho_1,\rho_2\in\sR$ with $\rho_1\ne\rho_2$.
We shall show that $(\rho_1,\rho_2)\in S$.
Let $0, x_1, x_2,\dots,x_m$ be the vertices, taken in
order, in the (finite) intersection of 
$\rho_1$ and $\rho_2$,
and write $x=x_m$, write
$y$ for the next vertex on $\rho_2$
as one moves from $x_m$ towards infinity, and
let $e=\langle x,y\rangle$. 

Let $\om\in\Om$ be the configuration which takes the value 1 on edges
$f\in \rho_1\cup\rho_2$, 
and the value 0 on all other edges.

Write the edges of $\rho_j^n$ in order as $f_j^n(1), 
f_j^n(2),\dots$, and those of $\rho_j$ as $f_j(1),
f_j(2),\allowmathbreak \dots$. 
By \refer{6.rholim},  for $j=1,2$,
for $i\ge 1$, and for all large $n$, we have that
$f_j^n(i) = f_j(i)$. Therefore, for all boxes $\La$,
and for $j=1,2$ and all large $n$,
the intersection of $E_\La$  with $\rho_j^n$ 
equals its intersection with $\rho_j$.
By the assumption that $(\rho_1^n,\rho_2^n)\in S$ 
for each $n$,
we have that $\om\in \K_{e,\La}^S$ for 
all $\La$, and therefore, by \refer{6.simcl},
$$
\om\in \lim_{\La\uparrow V}\K_{e,\La}^S = K_e^S.
$$
Since the only open rays in $\om$ 
are the sub-rays of $\rho_1$ and $\rho_2$, 
this can happen only if $(\rho_1,\rho_2)\in S$.
Therefore $S$ is closed,
as required.
\qed\enddemo

We return now to the universe of {\it equivalence\/}
relations on $\sR$, which we think of as binary
relations and denote by $\sim$.
We call ${\sim}\in\sE$ {\it closed\/}
(\resp, {\it open\/})
if it is closed (\resp, open) when viewed as a relation. The set
of closed  (\resp, open) 
equivalence relations is denoted $\sEc$  (\resp, $\sEo$).
By definition, $\sEc \subseteq \sSRc$ and $\sEo\subseteq \sSRf$.  Let $\sEcu$ 
($\subseteq
\sEf$) 
be the set of equivalence relations which,
when viewed as relations, are cutset relations.
The word `measurable' applied to an equivalence relation has
been explained in Section 4 and elaborated in Theorem 5.1.
We continue with a theorem concerning open equivalence relations.

\proclaim{Theorem \ch.5}
Let $\sim$ be an equivalence relation on $\sR$, and let 
$\sR(k)$ be the set of all paths of length $k$ 
starting at $0$.
The following are equivalent.
\ii{\rm(a)} $\sim$ is open.
\ii{\rm(b)} Each equivalence class of $\sim$ is an open 
subset of $\sR$.
\ii{\rm(c)} There exists an integer $k\ge 1$, and an
equivalence relation $\sim_k$ on $\sR(k)$, such that
$\sim$ is specified by $\sim_k$ in the sense that{\rm:}
$$
\rho_1\sim \rho_2 \q\text{if and only if}\q
\rho_1(k) \sim_k \rho_2(k),
$$
where $\rho(k)$ denotes the path comprising
the first $k$ edges of an 
infinite path $\rho$ from $0$.
\endproclaim

\demo{Proof}
Suppose that ${\sim}\in\sEo$. Then the sections 
$A_\rho=\{\rho'\in\sR: \rho\sim\rho'\}$,
$\rho\in\sR$,  are open,
and thus the equivalence classes are open.

Suppose next that (b) holds. The equivalence
classes cover the compact space $\sR$, whence 
there exists a finite sub-cover. Since the equivalence classes are
disjoint, no proper subset covers $\sR$, 
and therefore there exist only finitely many 
equivalence classes. Each class is the complement of the 
union of the others, and is therefore closed and thus 
compact.  Being open, each class is a union
of sets  of the form $\sR_v'$, $v\in V$
(since these form a base), and, being
compact, is a finite union of such sets. 
There exists therefore an integer $k$ such that every
equivalence class is a union of sets of the form
$\sR_w'$ as $w$ ranges over the set of
vertices at distance $k$ from the root 0.

That (c) implies (a) is obvious.
\qed\enddemo

\proclaim{Corollary \ch.6}
Every open equivalence relation is closed, and thus measurable. 
\endproclaim

\demo{Proof}
Let ${\sim}\in\sEo$. By Theorem \ch.5, the set
$\{(\rho,\rho')\in \sR^2: \rho\sim\rho'\}$ is a finite
union of closed sets and is therefore closed.
\qed\enddemo

We finish this section with a note.
There are or course many equivalence relations which are not closed. However,
with each relation $\sim$ may be 
associated a closed relation, termed the
{\it closure\/} of $\sim$ and denoted $\ol\sim$. We define $\ol\sim$ to be the
intersection of all closed equivalence relations $\sim'$
satisfying ${\sim}\le{ \sim'}$ with respect to the partial order
of 
\refer{4.partord}. It is easily checked that
$\ol\sim$ is itself a 
closed equivalence relation.

\endsection

\def\ch{6}
\section{\ch. Random-cluster measures with general 
boundary conditions}
Let us consider the tree $T=T_m'$ with
$m \ge 2$. 
We show first that 
$\sR^\sim_{p,q}\ne \es$ 
for ${\sim}\in\sEc$ and $q\ge 1$.
The proofs are given later in the section, and do not appear to
extend to the case $q<1$.

\proclaim{Theorem \ch.1} Let ${\sim}\in\sEc$
and let $0\le p\le 1$.
\ii{\rm(a)} Let $\pi$ be given by \refer{3.pprime},
and suppose $p \ne 1$, $q \ne 1$. 
The product measure $\phi_\pi$ satisfies 
$\phi_{\pi}
\in\sR^{\sim}_{p,q}$ if and only if 
$$
\phi_\pi(\text{\rm there exist two 
or more equivalent open rays})=0.
$$
\ii{\rm(b)} If $q\ge 1$, the weak limit
$$
\phi^{1,\sim}_{p,q}=\lim_{\La\uparrow V} \phi_{\Lpq}^{1,\sim}
$$
exists and satisfies 
$\phi^{1,\sim}_{p,q}\in\sR^{\sim}_{p,q}$.
\ii{\rm(c)} Let $q \ge 1$ and let 
$\phi\in\sR^{\sim}_{p,q}$.
Then
$$
\phi_\pi \lest \phi \lest \phi^{1,\sim}_{p,q}\lest \phi^{1}_{p,q}.
$$
\ii{\rm(d)} Let $q \ge 1$ and let 
$p$ be such that $\t^1(p,q)=0$; see \refer{4.I}.
The set $\sR^{\sim}_{p,q}$
comprises the singleton measure $\phi_{\pi}$ only.
\endproclaim

More generally, the proof of (c) shows that
$\phi^{1,\sim}_{p,q}\lest \phi^{1,\smash{\sim'}}_{p,q}$
whenever ${\sim}\le{\sim'}$. 
Part (b) may be extended as follows to arbitrary 
equivalence relations. 

\proclaim{Theorem \ch.2}
Let  ${\sim}\in \sE$, and let $\ol\sim$
be the closure of $\sim$. We have when $q\ge 1$ that
$$
\phi^{1,\sim}_{\La,p,q}\Rightarrow \phi_{p,q}^{1,\ol\sim}
\qq\text{as } \La\uparrow V.
\TAG
$$\label{6.genlim}
\endproclaim

This leads to the question: for ${\sim}\in\sEm$, 
is $\phi_{p,q}^{1,\ol\sim}$
a $(\sim)$-\rc\ measure?  The answer can be positive or negative,
as illustrated by the following examples.
\ii{1.} Let $\sim'$ be a closed relation with some 
equivalence class $C$ satisfying $|C|\ge 2$,
and let $\rho\in C$. Let $\sim$ be the relation
having the same equivalence classes as $\sim'$ except that $C$ is
replaced by the two equivalence classes $C\sm\{\rho\}$ and 
$\{\rho\}$.
It is easily seen that $\sim$ and $\sim'$ ($=\ol\sim$)
generate the same
family of \rc\ measures, and hence $\phi_{p,q}^{1,\sim}
\in\sR_{p,q}^\sim$ by \refer{6.genlim}.
(We return to such constructions after
Theorem 6.4.)
\ii{2.} Consider for definiteness the tree $T_2'$. Each of the three
sub-trees with root 0 is a binary tree, and each vertex therein
may be viewed as a combination of leftwards and rightwards steps from 0.
Let $\sT$ be the set of all (infinite) 0-rays which have only finitely
many leftwards steps. The set $\sT$ is dense in $\sR$, and
is countable and hence Borel. 
Let $\sim$ be the equivalence relation
whose unique non-trivial equivalence class is $\sT$. It is
easily seen that ${\ol\sim}={\sim^1}$. Let $q\ge 1$,
$p<1$, and $\phi\in\sR^{\sim}_{p,q}$.
By stochastic domination,  $\phi\lest \phi_p$, and therefore,
since $\sT$ is countable,
$$
\phi(\text{some ray in $\sT$ is open})=0.
$$
This implies that, with $\phi$-probability 1, there is no
pair of distinct equivalent open rays in the tree $T_2'$.  
By the forthcoming Lemma 6.5(a), the unique $(\sim)$-\rc\ 
measure is the product
measure $\phi_\pi$. By \refer{6.genlim}, 
$$
\phi_{p,q}^{1,\ol\sim} = \phi^1_{p,q}
= \lim_{\La\uparrow V}\phi_{\La,p,q}^{1,\sim} \ne \phi_\pi
$$ 
if
the conditions of \refer{4.I} hold. Therefore, $\phi_{p,q}^{1,\ol\sim}\notin
\sR_{p,q}^\sim$ under these conditions.

We turn next to the question of the uniqueness of \rc\ measures
for large $p$. 
We shall prove that $\sR_{p,q}^\sim$
is a singleton when ${\sim}\in\sEf$, $q\ge1$, 
and $p$
is sufficiently large. This extends [\ci{Jon99}, Thm 1.3],
which was concerned with
 the wired (maximal) boundary condition.
We leave as a problem the question whether this result extends to
arbitrary closed relations.

\proclaim{Theorem \ch.3} 
Let ${\sim}\in\sEf$ and let $0\le p\le 1$, $q \ge 1$. 
If $\pi\ge \pg(m)$ where 
$\pi$ and $\pg(m)$ are given in \refer{3.pprime}
and \refer{2.exactpg},
then $\sR^\sim_{p,q}$ comprises the singleton
$\phi^{1,\sim}_{p,q}$ only.
\endproclaim

Finally we turn to the question of the degree to which
the measure
$\phi_{p,q}^{1,\sim}$ characterises the relation $\sim$,
and for simplicity we consider first the case when 
$\sim$ is a cutset relation. 
Suppose that
${\sim}\in\sEcu$, which is to say that
${\sim}={\sim}^W$ for some cutset $W$. 
Let $T_W$ be the set of edges $e=\langle u,v\rangle$
which belong to no ray in 
$\bigcup_{w\in W}\sR_w$.
It is an easy consequence of the forthcoming Lemma \ch.5 that any
member $\phi$ of 
$\sR_{p,q}^\sim$ may be written in the form
$$
\phi=
\biggl\{\prod_{e\in T_W} \phi_{e,\pi}\biggr\}
\times\biggl\{\prod_{w\in W} \phi_{w,p,q}\biggr\},
\TAG
$$\label{5.prod2}
where $\phi_{e,\pi}$ 
is the Bernoulli measure on
$\{0,1\}$ associated with the state of $e$,
and $\phi_{w,p,q}$ is a wired \rc\ measure on the 
graph induced by $\sR_w$ with root $w$; 
conversely, any such $\phi$, with any choice of wired \rc\ measures
$\phi_{w,p,q}$, belongs to $\sR_{p,q}^\sim$.
In other words, $\phi\in\sR^{\sim}_{p,q}$ may 
be described as product measure
on $T_W$ with density $\pi$, combined with a wired \rc\ measure on
the graph formed by the rays in each given $\sR_w$, 
$w\in W$. 
In particular,
$$
\phi_{p,q}^{1,\sim} =
\biggl\{\prod_{e\in T_W} \phi_{e,\pi}\biggr\}
\times\biggl\{\prod_{w\in W} \phi_{w,p,q}^1\biggr\},
\TAG
$$\label{5.prod3}
where $\phi_{w,p,q}^1$ is the maximal
wired \rc\ measure on
$\sR_w$.

We say that $\phi_{w,p,q}^1$ {\it possesses a product
component\/} if there exists a non-empty subset $F$ of the set $E(\sR_w)$
of edges of $\sR_w$ such that
$$
\phi_{w,p,q}^{1} =
\biggl\{\prod_{f\in F} \phi_{f,\pi}\biggr\}
\times \psi,
\TAG
$$\label{5.prod4}
for some probability measure $\psi$ on the set
of configurations of $E(\sR_w)\sm F$.
When $q>1$ and $\t^1(p,q)>0$ then, by the results
of [\ci{Hag96}], 
or as a consequence of Lemma \ch.5(b) below, 
$\phi_{w,p,q}^1$
possesses no product component. It follows by
\refer{5.prod3} that, for any cutset relation
$\sim^W$, the measure $\phi^{1,\sim}_{p,q}$
is characterised by the set $W$. That is,
$$
\phi^{1,\sim}_{p,q} \ne \phi^{1,\sim'}_{p,q}
\qq\text{for } {\sim, \sim'} \in \sEcu,\ {\sim}\ne {\sim'}.
$$
whenever $q>1$ and \refer{4.I} holds.

Note in passing that $\sR^\sim_{p,q}$ contains
a continuum of distinct measures whenever 
there is a continuum of distinct wired measures on the graph induced by 
any given $\sR_w$. See [\ci{Hag96}, Section 5].

The above conclusion is extended to open equivalence
relations in the next theorem.

\proclaim{Theorem \ch.4}
Let $q> 1$ and let $p$ be such that $p \ne 1$
and $\t^1(p,q)>0$,
see \refer{4.I}. 
For distinct members $\sim,\sim'$ of $\sEf$,
we have that
$\phi_{p,q}^{1,\sim} \ne \phi_{p,q}^{1,\sim'}$.
\endproclaim

The conclusion of this theorem
is false with $\sEf$ replaced by
$\sEc$. For example, let $W$ be a cutset and let $w\in W$.
Consider the closed (equivalence) relation $S$ given by
$$
S= \biggl\{\bigcup _{x\in W,\ x\ne w} \sR_x' \times \sR_x'\biggr\}\cup\bigl\{(\rho,\rho): \rho\in\sR_w'\bigr\}.
$$
Now let $\rho,\rho'\in \sR_w'$ be distinct, and  let
$S' = S\cup\{(\rho,\rho'),(\rho',\rho)\}$. Then $S$ and $S'$ are
closed equivalence relations
which generate the same family of \rc\ measures.

The proofs of the three theorems above
are preceded by a lemma, a special case 
of which may be found in [\ci{Hag96}],
see also [\ci{GHM}]. We let 
$J_e=\{\text{$e$ is open}\}$. The $\s$-field
generated by the states of edges $f\in E$ 
with $f\ne e$ is denoted $\sT_e$.
The event $K_e^\sim$ is defined prior to Theorem 5.1.

\proclaim{Lemma \ch.5}
Let ${\sim}\in\sEm$ and $0\le p\le 1$, $q>0$.
\ii{\rm(a)} A probability measure $\phi$ on $(\Om,\sF)$
satisfies $\phi\in \sR^{\sim}_{p,q}$ if and
only if, for all $e\in E$ and for $\phi$-almost
every $\xi$,
$$
\phi(J_e\mid \sT_e)(\xi)
=\cases p &\text{if } \xi\in K_e^\sim,\\
\pi &\text{if } \xi\notin K_e^\sim,
\endcases
\TAG
$$\label{5.sngle}
where $\pi$ is given in \refer{3.pprime}
\ii{\rm(b)} Let $e\in E$ and $\phi\in \sR^\sim_{p,q}$ where 
$p\neq0,1$ and
$q\ne 1$.
The (random) state of $e$ is independent of the states of
$E\sm\{e\}$,  equalling $1$ with probability $\pi$, if
and only if
$$
\phi(K_e^\sim) =0.
$$
\endproclaim

\demo{Proof}
(a) By an application of 
[\ci{Lig}, Propn IV.1.8],
the \rc\ measure $\phi_{G,p,q}$ on a finite graph $G=(V,E)$
is characterised by the statement that, for all edges
$e=\langle x,y\rangle\in E$,
$$
\phi_{G,p,q}(J_e\mid \sT_e)(\xi) 
=\cases p &\text{if } \xi\in K_e,\\
\pi &\text{if } \xi\notin K_e,
\endcases
\TAG
$$\label{5.sngle2}
where $K_e$ is the event that $x$ and $y$ are joined
by an open path of $E\sm\{e\}$. See [\ci{G04}, Thm 2.1], for
example. 
If $\phi$ satisfies the condition of the lemma then, by Definition
\refer{3.defrc}, $\phi\in \sR^\sim_{p,q}$. 
The converse is similar.

\flushpar
(b)  When $p\ne0,1$ and $q \ne 1$, we have that $p\ne \pi$, and the
claim follows by part (a).
\qed\enddemo

\demo{Proof of Theorem \ch.1}
(a)
Suppose $p\ne 0$; the case $p=0$ is trivial.
If $\phi_\pi$ satisfies the condition, then
$\phi_\pi(K_e^\sim) = 0$ for all $e\in E$, implying
that $\phi_\pi$ satisfies \refer{5.sngle}.
By Lemma 6.5(a),
$\phi_\pi\in\sR^\sim_{p,q}$. The converse argument
is valid when $p\ne \pi$, which requires that
$p\ne 1$, $q\ne 1$.

\flushpar(b) Let ${\sim}\in\sEc$.
By the FKG inequality in the usual way
(see [\ci{G93}, Thm 3.1(a)] for example), 
the limit 
$$
\phi=\lim_{\La\uparrow V} \phi_{\Lpq}^{1,\sim}
$$
exists and is a probability measure.
We prove next that $\phi$ satisfies
\refer{5.sngle} for all $e$ and $\phi$-almost
every $\xi$, and the claim will follow.

For $\xi\in\Om$ and $F\subseteq E$, write $[\xi]_F$
for the set of all configurations which agree 
with $\xi$ on $F$. For $W\subseteq V$ and $e\in E_W$, let
$[\xi]_W$ (\resp, $[\xi]_{W\sm e}$) be an abbreviation for 
$[\xi]_{E_W}$ (\resp, $[\xi]_{E_W\sm\{e\}})$.
We write $\xi_W^1$ for the configuration
which agrees with $\xi$ on $E_W$ and which equals 1 elsewhere.
For economy of notation,
we shall omit explicit reference to the values
of $p$ and $q$ in the rest of this proof, and thus
we write $\phi_\La^1 = \phi_{\La,p,q}^{1,\sim}$.

By the martingale
convergence theorem, for $e=\langle x,y\rangle
\in E$ and
$\phi$-almost every $\xi$,
$$
\align
\phi(J_e\mid \sT_e)(\xi) &= \lim_{\La\uparrow V} 
 \frac {\phi(J_e,\ [\xi]_{\La\sm e})}
 {\phi([\xi]_{\La\sm e})}\TAG\\\label{5.big}
&=\lim_{\La\uparrow V} \lim_{\De\uparrow V}
 \frac {\phi_\De^1(J_e,\ [\xi]_{\La\sm e})}
 {\phi_\De^1([\xi]_{\La\sm e})}\\
&=\lim_{\La\uparrow V} \lim_{\De\uparrow V}
\phi_\De^1\bigl(\phi_\De^1(J_e\mid [\xi]_{\De\sm e})
\bigmid [\xi]_{\La\sm e}\bigr)\\
&=\lim_{\La\uparrow V} \lim_{\De\uparrow V}
 \phi_\De^1(g_\De\mid [\xi]_{\La\sm{e}}),
\endalign
$$
by \refer{5.sngle2}, where
$$
g_\De(\xi) = \pi + (p-\pi)1_{\K^\sim_{e,\De}}(\xi)
$$
and $\K^\sim_{e,\De}=\{\om\in\Om:
\om_\De^1\in K_e^\sim\}$. Here and later, $1_A$
denotes the indicator function of the event $A$.
Since $\sim$ is closed
by assumption, we have by Theorem 5.4 that
$$
\K^\sim_{e,\De} \downarrow K_e^\sim\qq\text{as } \De \uparrow V,
\TAG
$$\label{5.condv}
whence $g_\De\downarrow g$ where $g = \pi+(p-\pi)1_{K_e^\sim}$.

We claim that
$$
\phi_\De^1(g_\De\mid [\xi]_{\La\sm{e}})
\to \phi(g\mid [\xi]_{\La\sm{e}})
\qq\text{as } \De\uparrow V,
\TAG
$$\label{5.big2}
and we prove this as follows. Let $\De'$ be a box satisfying
$\La\subseteq \De'\subseteq \De$,
and write $\psi_\De(\cdot)=\phi_\De^1(\cdot\mid [\xi]_{\La\sm{e}})$.
Since $g_\De$ is non-increasing in $\De$, 
$$
\psi_{\De}(g) \le \psi_{\De}(g_\De)\le \psi_{\De}(g_{\De'}).
$$
We take the limits as $\De\uparrow V$ and
$\De'\uparrow V$, in that order, and we
appeal to the
dominated convergence theorem to deduce \refer{5.big2}.

By the martingale convergence theorem,
$$
\phi(g\mid [\xi]_{\La\sm{e}})\to g(\xi)\qq
\text{as $\La\uparrow V$, for $\phi$-almost every $\xi$}.
$$ 
Hence, \refer{5.big} and \refer{5.big2} show that \refer{5.sngle}
holds, and the result follows.

\flushpar(c) 
Let $A$ be an increasing cylinder
event in $\sF$, and suppose that $A$ is defined on a finite set $B$
of edges.
Let $\La$ be a box  of the tree satisfying $E_\La\supseteq B$.
In the construction of 
$\phi_{\Lpq}^{\xi,\sim}$ in Section 4, we add a certain set of (permanently
open) new edges to $\La$; for 
$\phi_{\Lpq}^{1,\sim}$ we add a larger (or equal) set of new edges and for 
$\phi_{\Lpq}^{1}$ an even larger set; on the other hand, no such edges
are added for 
$\phi_{\Lpq}^{0}=\phi_{\Lpq}^{0,\sim_0}$.
By the FKG inequality, the addition of a new open edge gives rise to a
stochastically larger \rc\ measure. Thus, for every $\xi$,
$$
\phi_{\Lpq}^{0}(A)
\le
\phi_{\Lpq}^{\xi,\sim}(A)
\le
\phi_{\Lpq}^{1,\sim}(A)
\le
\phi_{\Lpq}^{1}(A).
$$
Here, $\phi_{\Lpq}^{0}$ is the product measure $\phi_\pi$, cf.\ Theorem
4.1(a).
Using \refer{3.rcdef} and taking the expectation over $\xi$ we find
$$
\phi_\pi(A) 
\le 
\phi(A) 
\le
\phi_{\Lpq}^{1,\sim}(A)
\le
\phi_{\Lpq}^{1}(A).
$$
Now let $\La\uparrow V$.

\flushpar(d) 
Let $\phi\in \sR^\sim_{p,q}$.
Since $\t^1(p,q)=0$, 
there is $\phi^1_{p,q}$-a.s.\ no open ray.
By (c), there is thus $\phi$-a.s.\ no open ray, and thus, by
\refer{3.defrc}, for every cylinder set $A$ defined on some finite edge-set
$B$ and every box $\La$ with $E_\La\supseteq B$,
$$
\phi_{\Lpq}^{\xi,\sim}(A)
=
\phi_{\Lpq}^{\xi,\sim_0}(A)
=
\phi_{\pi}(A)
\qquad
\text{for $\phi$-a.e.\  $\xi$}.
$$
Hence, by \refer{3.rcdef} and taking the expectation, $\phi(A)=\phi_\pi(A)$.
\qed\enddemo

\demo{Proof of Theorem 6.2}
Let ${\sim}\in\sE$. Using \refer{4.stord} in the usual way,
we may restrict ourselves for simplicity to sets $\La$ which are boxes.
Let $\La$ be a finite box of $V$.
For $v\in\La$, let $C_v$ be the subset of $\sR$ containing
all 0-rays whose final point of intersection with $\La$ is $v$.
Let $\approx$ be the lowest equivalence relation (in
the sense of the natural partial order on
equivalence relations, see \refer{4.partord}) on the $C_v$ 
such that:
$C_v \approx C_w$ if there exist $\rho_v\in C_v$,
$\rho_w\in C_w$ such that $\rho_v\sim \rho_w$.
The relation $\approx$ may be extended without change of notation
to an equivalence relation on $\sR$.
By construction, $\approx$ is closed (since
each $C_v$ is closed). Since
${\sim}\le {\approx}$, we have that 
$$
{\sim}\le{\ol\sim}\le{\approx}.
\TAG
$$\label{6.equord}

Let $\om\in\Om$, and let $k^{1,\sim}(\La,\om)$ be given as before
\refer{3.defrc}. We claim that
$$
k^{1,\sim}(\La,\om) = k^{1,\approx}(\La,\om).
\TAG
$$\label{6.suff}
By \refer{6.equord}--\refer{6.suff}, 
$k^{1,\sim}(\La,\om)=k^{1,\ol\sim}(\La,\om)$, and hence
$\phi_{\La,p,q}^{1,\sim} = 
\phi_{\La,p,q}^{1,\ol\sim}$. The claim of the theorem follows
on taking the limit $\La\uparrow V$.

We prove \refer{6.suff} next. Let $A^\sim= A^\sim(\La,\om)$
be the set of edges which are added to $G(\La,\om)$
in the construction
prior to \refer{3.defrc} of the measure $\phi^{1,\sim}_{\La,p,q}$. 
Since ${\sim}\le{\approx}$, $A^\sim \subseteq A^\approx$.
It therefore suffices to prove that, if $\la u,v\ra\in A^\approx$,
there exists a path of $A^\sim$ joining $u$ to $v$.
Suppose $\la u,v\ra \in A^\approx$, so that
$C_u \approx C_v$. By the definition
of $\approx$, there exist vertices $u_0=u, u_1,
\dots, u_k=v$ such that, for $0\le i < k$,
$C_{u_i}$ and $C_{u_{i+1}}$ contain two
$\sim$-equivalent rays. 
Therefore, $\la u_i,u_{i+1}\ra\in A^\sim$, and
hence $u$ and $v$ are connected by a path
in $A^\sim$. The proof is complete.
\qed\enddemo

\demo{Proof of Theorem \ch.3}
This is inspired by the proof of 
[\ci{G93}, Thm 5.3(c)], see also [\ci{Jon99}, Thm 1.3].
For any cutset $C$, we write $\ins(C)$ for the set of all vertices
reachable from 0 along paths of $T$ disjoint from $C$,
and $\ol C= C\cup\ins(C)$. The sub-$\s$-field 
of $\sF$ generated by the states of
edges having no more than one vertex in $\ol C$ is 
denoted $\sT_C$.

Let $\pi\ge\pg(m)$ and ${\sim}\in\sEf$. 
Choose $k$ according
to Theorem 5.5(c), and let $W$ be the set of
all vertices at distance $k$ from 0.
Let $A$ be an increasing cylinder
event in $\sF$, and suppose that $A$ is defined on the set of 
edges in some box $B$ of the tree which we may take 
sufficiently large to contain $W$.  Let 
$\phi\in\sR^{\sim}_{p,q}$. We shall show that
$$
\phi(A) = \phi^{1,\sim}_{p,q}(A),
\TAG
$$\label{5.enuf}
and the claim will follow.

Let $\La$ be a box satisfying $\La\supseteq B$ and
let $\om\in\Om$. 
Recall from Section 2 that a vertex $x$ is blue if there exists
an $x$-ray that is open (for the configuration $\om$). 
For any box 
$\De\supseteq\La$, let $\sB_\De= \sB_\De(\om)$ be the 
set of blue cutsets contained in $\De$, noting that $\sB_\De=\es$
is possible. 
There is a natural partial order
on $\sB_\De$ given by $C_1\le C_2$ if $C_1\subseteq
\ol{C_2}$. Let $\cmax$ be the maximal blue cutset in this partial ordering,
thus,
$$
\ol{\cmax} = \bigcup_{C\in \sB_\De} \ol C.
$$
Note that, for any cutset
$C \subseteq \De$,
the event 
$\edc=\{\cmax = C\}$ lies in the $\s$-field
$\sT_C$.

Let $\eps>0$. There exists a deterministic box 
$\De'=\De'(\La,\eps)\supseteq \La$ such that:
$$
\phi_{\pi}(E_\De)
>1-\eps
\q\text{for all $\De\supseteq \De'$},
\TAG
$$\label{5.prod}
where $E_\De$ is the event
$$
E_\De = \{\text{$\sB_\De\ne\es$ and  
$\La\subseteq \ins(\cmax)$}\}
=\bigcup_{\ins(C)\supseteq\La} \edc
.
$$
Corollary 2.3 and the assumption $\pi\ge \pg(m)$ have
 been used here. The corollary was phrased for the tree $T_m$
rather than $T_m'$, but it is easily seen to be
valid for either tree.

Let $C$ be a cutset with $\La\subseteq\ins(C)$, $C\subseteq\De$.
On the event $\edc$,
$C$ is a blue cutset. By \refer{3.rcdef} and the fact that
$W\subseteq B\subseteq \ins(C)$,
$$
\phi(A\mid \sT_C) = \phi_{\ol C,p,q}^{1,\sim}(A) \qq
\text{$\phi$-a.s., on $\edc$}.
$$
By the FKG inequality (see \refer{4.stord}), 
$$
\phi_{\De,p,q}^{1,\sim}(A) \le \phi(A\mid \sT_C) 
\le \phi_{\La,p,q}^{1,\sim}(A)\qq
\text{$\phi$-a.s., on $\edc$}.
$$
We take expectations and sum over cutsets $C\subseteq \De$ such that
$\La\subseteq\ins(C)$    
to find that 
$$
\phi_{\De,p,q}^{1,\sim}(A)\phi(E_\De) 
 \le \phi(A) -\phi(A,\, \text{not } E_\De)
\le \phi_{\La,p,q}^{1,\sim}(A)\phi(E_\De).
\TAG
$$\label{6.*}
By \refer{5.prod} and the fact 
from Theorem \ch.1(c) 
that 
$\phi\gest \phi_\pi$,
if $\De$ is large enough,
$$
0\le 
\phi(A,\, \text{not } E_\De)
\le
1-\phi(E_\De) 
\le 1-\phi_\pi(E_\De)
 < \eps.
$$
We pass to the limits in \refer{6.*} as $\De\uparrow V$, 
$\eps \downarrow 0$,   
and $\La\uparrow V$
to obtain \refer{5.enuf} as required.
\qed\enddemo

\demo{Proof of Theorem \ch.4}
Let $|x|$ denote the length of the unique path of $T_m'$
from 0 to $x$. Fix $k\ge 0$, let $\La_k$ be the set of all $x$
with $|x|\le k$, and let $\pd\La_k=\La_k\sm\La_{k-1}$ as usual.
We call $y$ a {\it descendant\/} of $x$ if
the unique path from 0 to $y$ passes through $x$.
For $x\in V$, write $D_m(x)$ for the set of all descendants $y$ of $x$
with $|y|=m$, and $D(x)=\bigcup_{m >|x|} D_m(x)$.

Choose ${\sim}\in\sEf$, and fix $k=k(\sim)$ such that the conclusion of
Theorem 5.5(c) holds. Let $\sim_k$ be given as in that theorem.
For $x,y\in\pd\La_k$, we
write $x\sim_k y$ if the unique paths
$\pi_{0,x}$ (\resp, $\pi_{0,y}$)
 from 0 to $x$ (\resp, $y$)
satisfy $\pi_{0,x}\sim_k\pi_{0,y}$.
We have that
$$
\phi_{p,q}^{1,\sim}= \lim_{n\to\oo} \psi^\sim_{n},
$$
where $\psi_n^\sim=\phi_{\La_n,p,q}^{1,\sim}$. Let $n  >k=k(\sim)$.
The measure $\psi_n^\sim$ 
may be considered as the \rc\ measure  
on the graph $G_n^\sim=(\La_n^\sim, E_{\La_n})$ where $\La_n^\sim$ is
obtained from $\La_n$ by identifying any pair of vertices
$u,v\in\pd\La_n$ having the property that there exist
$\rho_u\in\sR_u$, $\rho_v\in\sR_v$ with
$\rho_u\sim\rho_v$. Since $n > k$, this implies two levels
of identifications:
\ii{(i)} for every $x\in\pd\La_k$, we identify the set $D_n(x)$
of vertices,
\ii{(ii)}  if $x,y\in\pd\La_k$ are such that $x\sim_k y$,
we identify all vertices in $D_n(x)\cup D_n(y)$.
\flushpar
Part (ii) incorporates part (i), since $x\sim_k x$ for $x\in\pd\La_k$,
but we express it thus in order to emphasize the role
of the equivalence relation $\sim$.

We recall a basic fact (see
[\ci{G04}], especially the appendix). Consider the \rc\ measure
on a finite graph $G=(A,B)$ with cluster-weighting factor $q$ and
a family $\bop=$ $ =(p_b: b\in B)$ of edge-parameters. 
Any set $B'$ of one or more edges in parallel
(\resp, in series) may be replaced by a single edge with an
associated edge-parameter that is a function of $(p_b: b \in B')$
and $q$. We call such a replacement a {\it parallel\/}
(\resp, {\it series\/}) replacement. 

{}From the finite sub-tree $(\La_k,E_{\La_k})$ of $T_m'$ we construct
as follows a further tree denoted $T$.
To each $x\in\pd\La_k$ we attach
a further
edge $[x,x'\rangle$ to the vertex $x$, arriving thus at a tree which
we denote $T$ and which is illustrated in Figure \refer{6.fig2}.
These new edges are called {\it attachments\/}. From $T$ we obtain
the graph $T^\sim$ by identifying $x'$ and $y'$ whenever $x\sim_k y$.
The graph $G_n^\sim$ may be transformed by a sequence of parallel and series
replacements into $T^\sim$. The \rc\ measure $\psi_n^\sim$ corresponds
to the \rc\ measure $\nu_n^\sim$ on $T^\sim$ for which the attachment edge 
$[x,x'\rangle$ has an associated parameter $p_n$ that depends
only on $p$, $q$, $k$, and $n$, and not further on the choice of 
$x$ and $\sim$.
Since $\psi_n^\sim\Rightarrow \phi_{p,q}^{1,\sim}$, we have that
$p_n\to p_\oo$ for some $p_\oo\in[0,1]$, and that $\phi_{p,q}^{1,\sim}$
corresponds to the \rc\ measure $\nu_\oo^\sim$ on $T^\sim$ for which
the attachment edges have associated parameter $p_\oo$.  

\topinsert
\figure
\mletter{e_x}{4.9}{3}
\mletter{x}{4.68}{3.7}
\mletter{0}{6.4}{2.1}
\lastletter{x'}{3.45}{3.9}
\centerline{\epsfxsize=5cm
  \epsfbox{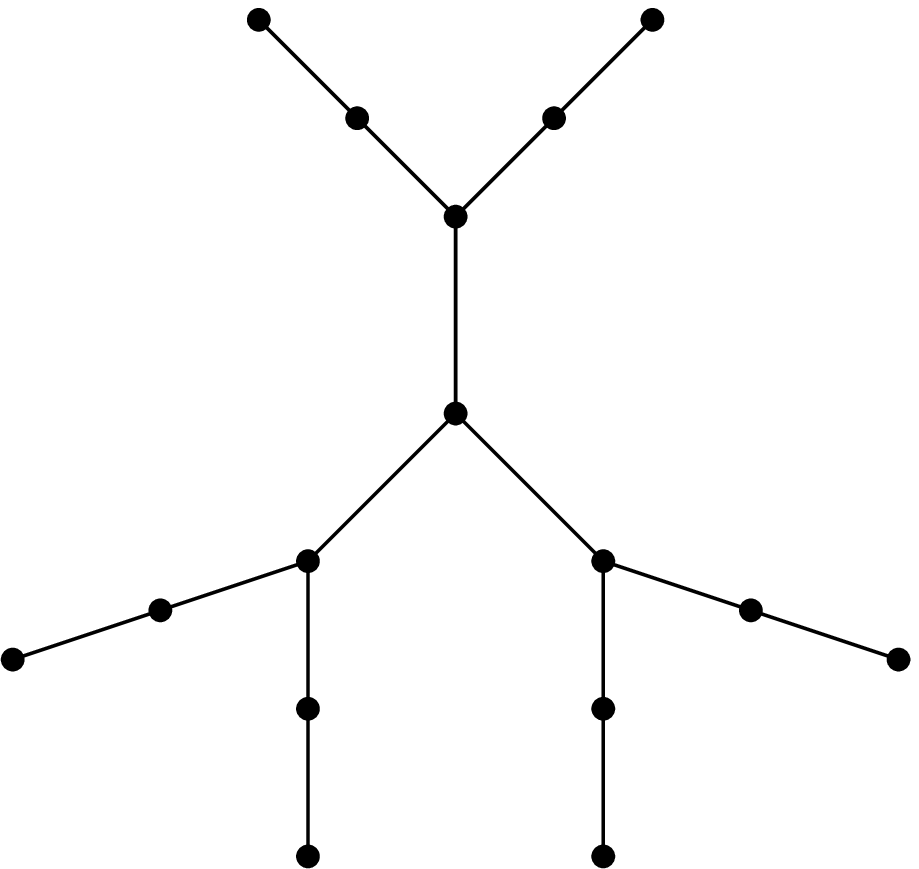}}

\capt{To each boundary vertex $x$ of the first $k$ 
($=2$) generations of $T_m'$
is attached a new edge $[x,x'\rangle$.}
\label{6.fig2}
\endinsert

Assume now that $p$ is such that $0<\t^1(p,q)<1$. It is easily deduced 
that $p_\oo \in(0,1)$. For $x\in\pd\La_k$, we write 
$e_x$ for the unique edge of $T_m'$ of the form
$[z,x\rangle$. The marginal (joint) law of the states of the
edges $(e_x: x\in\pd\La_k)$ is the same under $\phi_{p,q}^{1,\sim}$
as under $\nu_\oo^\sim$. 

Let $q>1$,
and let $x,y\in\pd\La_k$, $x\ne y$. 
Let $C$ be the event that $e_z$ is closed for all $z\in\pd\La_k\sm\{x,y\}$.
The conditional measure given $C$ of $\nu_\oo^\sim$ is a \rc\ measure
on the graph $T^\sim$ with the edges $e_z$, $z\in\pd\La_k\sm\{x,y\}$,
removed, and we denote this graph by $T^\sim_{x,y}$. 
It is a fundamental 
property of \rc\ measures with $q\ne 1$ on a finite graph,
with edge-parameters lying in $(0,1)$,
that the states of two edges are dependent random variables
if and only if there exists some cycle containing both edges
(see [\ci{G04}]).
There exists a cycle in $T^\sim_{x,y}$ containing both $e_x$
and $e_y$ if and only if $x\sim_k y$.
Therefore, $x\sim_k y$ if and only if
$$
\nu_\oo^\sim(J_x \mid  J_y\cap C)
\ne
\nu_\oo^\sim(J_x \mid  J_y^{\text{c}}\cap C),
\TAG
$$\label{6.last}
where $J_x=J_{e_x}$ is the event that $e_x$ is open.

The proof is nearly complete. Let $\sim$ and $\sim'$ be distinct
open equivalence relations, and choose $k$ sufficiently large that
the conclusion of Theorem 5.5(c) holds for both $\sim$ and $\sim'$.
Since $\sim$ and $\sim'$ are distinct, there exist $x,y\in\pd\La_k$
such that $x\sim_k y$ but $x\nsim'_k y$. Therefore,
$\nu_\oo^\sim$ satisfies \refer{6.last} but $\nu_\oo^{\sim'}$
does not, and hence $\phi_{p,q}^{1,\sim} \ne \phi_{p,q}^{1,\sim'}$
as required.
\qed\enddemo

\endsection

\section{Acknowledgements} 
This work was done during a visit by the second author
to Cambridge University with support from
the Royal Swedish Academy of Sciences,
the London Mathematical Society, and Churchill College, Cambridge.
\endsection

\enddocument